\documentclass[12pt,draft]{amsart}

\usepackage[cp1251]{inputenc}
\usepackage[T2A]{fontenc}
\usepackage[english,russian]{babel}
\usepackage{amsmath,amsfonts,amssymb}
\usepackage{geometry}
\author[В. А. Михайлец, А. А. Мурач]
{В. А. Михайлец (V. A. Mikhailets), А. А. Мурач (A. A. Murach)}

\title[Эллиптическая задача
в уточненной шкале пространств]{Регулярная эллиптическая краевая
задача в двусторонней уточненной шкале пространств \\
\vspace{0.5cm} Regular elliptic boundary-value problem in \\ a
two-sided refined scale of spaces}

\address{Институт математики НАН Украины,
ул. Терещенкивська 3, Kиев, Украина, 01601}

\email{mikhailets@imath.kiev.ua}

\address{Институт математики НАН Украины,
ул. Терещенкивська 3, Kиев, Украина, 01601 \\ \indent Черниговский
государственный технологический Университет, ул. Шевченка 95,
Чернигов 14027, Украина}

\email{murach@imath.kiev.ua}

\subjclass[2000]{35J40, 46E35}

\date{10/12/2007}

\keywords{Regular elliptic boundary-value problem, generalized
solution, the Fredholm property, local smoothness of a solution,
the H\"ormander spaces, regularly varying function, interpolation
with function parameter.}

\begin{document}

\maketitle

\begin{abstract}
A regular elliptic boundary-value problem over a bounded domain
with a smooth boundary is studied. We prove that the operator of
this problem is a Fredholm one in the two-sided refined scale of
the functional Hilbert spaces and generates a complete collection
of isomorphisms. Elements of this scale are the isotropic spaces
of H\"{o}rmander--Volevich--Paneah and some its modifications.
A~priori estimate for the solution is established and its
regularity is investigated.
\end{abstract}

\vspace{1cm}

\noindent\textbf{Введение.} В работах Ж.-Л. Лионса, Э. Мадженеса
\cite{LM71} и Ю. М. Березанского, С. Г. Крейна, Я. А. Ройтберга
\cite{BeKrRo, Ro64, Ber65, Ro96} установлены теоремы о полном
наборе изоморфизмов, кототорый осуществляет оператор регулярной
эллиптической краевой задачи в двусторонней шкале пространств
функций/распределений. Полнота набора означает, что указанные
изоморфизмы выполняются между пространствами
функций/распределений, которые имеют соответственно $s$ и $s-2q$
производных, где $s$ --- \textit{произвольное} вещественное число,
а $2q$ --- порядок оператора. Позитивная часть двусторонней шкалы
($s\geq2q$) состоит из пространств Соболева, а негативная часть
($s<2q$) --- из модифицированных специальным образом соболевских
пространств. Известны две такие модификации. Модификация
Лионса-Мадженеса состоит из некоторых \textit{сужений} соболевских
пространств, что позволяет обеспечить непрерывность краевых
операторов в них. Иной подход предложен Я.~А.~Ройтбергом
\cite{Ro64, Ber65, Ro96} и основывается на \textit{расширении}
соболевских пространств. Более того, обобщенное решение краевой
задачи трактуется как вектор, компоненты которого принадлежат
соболевским пространствам и связаны определенным образом между
собой. Это позволило исследовать эллиптическую краевую задачу,
правые части которой являются произвольными распределениями.

Теоремы о полном наборе изоморфизмов были доказаны Я. А.
Ройтбергом \cite{Ro96} также для нерегулярных эллиптических
краевых задач и краевых задач для эллиптических систем
дифференциальных уравнений. В наиболее общем виде они установлены
А.~Н.~Кожевниковым \cite{Kozh01} для псевдодифференциальных
эллиптических краевых задач. Эти теоремы имеют ряд важных
приложений (см. \cite{Ro96} и приведенную там литературу). Среди
них особое место занимают утверждения о повышении локальной
гладкости решения эллиптической краевой задачи. В этой связи
является актуальным изучение эллиптических задач в двусторонних
шкалах пространств, дающих более тонкую градацию гладкостных
свойств распределений, чем соболевская шкала. К их числу относятся
гильбертова шкала специальных изотропных пространств Хермандера --
Волевича -- Панеяха \cite{Her65, Her86, VoPa65, Pa00}
$$
H^{s,\varphi}:=H_{2}^{\langle\cdot\rangle^{s}\,\varphi(\langle\cdot\rangle)},
\quad\langle\xi\rangle:=\bigl(1+|\xi|^{2}\bigr)^{1/2},
$$
где $s\in\mathbb{R}$, а функциональный параметр $\varphi$ является
медленно меняющейся на $+\infty$ по Карамата функцией $t\gg1$. В
частности, допустима любая эталонная функция
$$
\varphi(t)=(\log t)^{r_{1}}(\log\log
t)^{r_{2}}\ldots(\log\ldots\log
t)^{r_{n}},\quad\{r_{1},r_{2},\ldots,r_{n}\}\subset\mathbb{R},\;n\in\mathbb{N}.
$$
Эта \textit{уточненная} шкала введена и изучена в \cite{MM2a,
MM2b}. Она содержит соболевскую шкалу
$\{H^{s}\}\equiv\{H^{s,1}\}$, привязана к ней числовым параметром
$s$, но много тоньше нее.

Пространства $H^{s,\varphi}$ очень естественно возникают в ряде
спектральных задач: сходимость спектральных разложений
самосопряженных эллиптических операторов почти всюду, по норме
пространства $L_{p}$ с $p>2$ или $C$ (см. обзор \cite{AIN76});
спектральная асимптотика общих самосопряженных эллиптических
операторов в ограниченной области, формула Г.\;Вейля, точная
оценка остаточного члена в ней (см. \cite{Mikh82, Mikh89}) и др.
Можно ожидать, что они окажутся полезными и в иных "тонких"\;
вопросах. Благодаря своим интерполяционным свойствам
$H^{s,\varphi}$ занимают особое место среди пространств обобщенной
гладкости, которые все активнее исследуются и используются в
последние годы (см. обзор \cite{KalLiz87}, недавние работы
\cite{HarMou04, FarLeo06} и приведенную там литературу).

В настоящей статье изучается регулярная эллиптическая краевая
задача в двусторонней уточненной шкале пространств, негативная
часть которой модифицирована по Ройтбергу. Доказано, что оператор
этой задачи ограничен, фредгольмов и порождает полный набор
изоморфизмов в такой шкале. Исследована уточненная локальная
гладкость решения эллиптической задачи. В качестве приложения дано
достаточное условие классичности обобщенного решения задачи.

Отметим для полноты изложения, что в позитивнной части уточненной
шкалы неоднородная эллиптическая краевая задача изучена ранее в
\cite{Sh74, MM2b, MM2c}. \textit{Полуоднородные} эллиптические
краевые задачи можно исследовать в \textit{двусторонних}
уточненных шкалах \textit{без} их модификации (см. \cite{MM5,
MM6}). Однако, неоднородная краевая задача в негативной части
шкалы не сводится к двум полуоднородным, т. к. их решения являются
распределениями разной природы. Случай эллиптических ПДО в
двусторонней уточненной шкале пространств на замкнутом компактном
могообразии исследован авторами в \cite{MM10, MM11, MM12}. Укажем
также на работы \cite{M94a, M94b}, где эллиптическая краевая
задача изучалась в двусторонних модифицированных шкалах
пространств Лизоркина-Трибеля и Никольского-Бесова.

\vspace{0.5cm}

\textbf{1. Постановка задачи и  основной результат.} Пусть
$\Omega$ --- ограниченная область в евклидовом пространстве
$\mathbb{R}^n$ $(n\geq2)$ с границей $\Gamma$, которая является
бесконечно гладким замкнутым многообразием размерности $n-1$.
Предполагается, что область $\Omega$ локально лежит по одну
сторону от $\Gamma$. Обозначим
$\overline{\Omega}=\Omega\cup\Gamma$.

Рассмотрим следующую неоднородную краевую задачу в области
$\Omega$:
$$
A\,u=f\;\;\text{в}\;\;\Omega,\quad
B_{j}\,u=g_{j}\;\;\text{на}\;\;\Gamma\;\;\text{при}\;\;
j=1,\ldots,q.\eqno(1.1)
$$
Здесь и далее $A$ --- линейное дифференциальное выражение в
$\overline{\Omega}$ произвольного четного порядка $2q\geq2$, а
$B_{j}$, где $j=1,\ldots,q$, --- граничное линейное
дифференциальное выражение на $\Gamma$ порядка $m_{j}\leq2q-1$.
Все коэффициенты выражений $A$ и $B_{j}$ являются
комплекснозначными функциями, бесконечно гладкими в
$\overline{\Omega}$ и на $\Gamma$ соответственно. Положим
$B:=(B_{1},\ldots,B_{q})$.

Всюду далее предполагается, что краевая задача (1.1)
\textit{регулярная эллиптическая}. Это означает \cite[с. 137 -
138]{LM71}, \cite[с. 167]{Kr72}, что выражение $A$ правильно
эллиптическое в $\overline{\Omega}$, а набор граничных выражений
$B$ нормальный и удовлетворяет условию дополнительности по
отношению к $A$ на $\Gamma$. Из условия нормальности вытекает, что
порядки $m_{j}$ граничных дифференциальных выражений все различны.

Наряду с задачей (1.1) рассмотрим краевую задачу
$$
A^{+}\,v=\omega\;\;\text{в}\;\;\Omega,\quad
B^{+}_{j}\,v=h_{j}\;\;\text{на}\;\;\Gamma\;\;\text{при}\;\;j=1,\ldots,q.\eqno(1.2)
$$
Она формально сопряжена к задаче (1.1) относительно формулы Грина:
$$
(Au,v)_{\Omega}+\sum_{j=1}^{q}\;(B_{j}u,\,C_{j}^{+}v)_{\Gamma}
=(u,A^{+}v)_{\Omega}+\sum_{j=1}^{q}\;(C_{j}u,\,B_{j}^{+}v)_{\Gamma},
\;\;u,v\in C^{\infty}(\,\overline{\Omega}\,).
$$
Здесь $A^{+}$ --- сопряженное к $A$ линейное дифференциальное
выражение порядка $2q$ с коэффициентами класса
$C^{\infty}(\,\overline{\Omega}\,)$, а $\{B^{+}_{j}\}$,
$\{C_{j}\}$, $\{C^{+}_{j}\}$ --- некоторые нормальные системы
линейных дифференциальных граничных выражений с коэффициентами
класса $C^{\infty}(\Gamma)$. Их порядки удовлетворяют условию:
$$\mathrm{ord}\,B_{j}+\mathrm{ord}\,C^{+}_{j}=
\mathrm{ord}\,C_{j}+\mathrm{ord}\,B^{+}_{j}=2q-1.
$$
Здесь через $(\cdot,\cdot)_{\Omega}$ и $(\cdot,\cdot)_{\Gamma}$
обозначены скалярные произведения в пространствах $L_{2}(\Omega)$
и $L_{2}(\Gamma)$ функций, интегрируемых с квадратом в $\Omega$ и
на $\Gamma$ соответственно, а также естественные расширения по
непрерывности этих скалярных произведений.

Положим
$$
N:=\{u\in C^{\infty}(\,\overline{\Omega}\,):\;Au=0\;\;\mbox{в}\;\;
\Omega,\;\;
B_{j}u=0\;\;\mbox{на}\;\;\Gamma\;\;\mbox{для}\;\;j=1,\ldots,q\},
$$
$$
N^{+}:=\{v\in
C^{\infty}(\,\overline{\Omega}\,):\;A^{+}v=0\;\;\mbox{в}\;\;
\Omega,\;\;
B^{+}_{j}v=0\;\;\mbox{на}\;\;\Gamma\;\;\mbox{для}\;\;j=1,\ldots,q\}.
$$
Поскольку задачи (1.1) и (1.2) являются регулярными
эллиптическими, пространства $N$ и $N^{+}$ конечномерны \cite[с.
191]{LM71}, \cite[с. 168]{Kr72}.

Для простоты формулировок предположим в этом пункте, что
$N=N^{+}=\{0\}$.

Напомним следующий классический результат \cite[с. 191]{LM71},
\cite[с. 169]{Kr72}: оператор $(A,B)$, соответствующий задаче
(1.1), определяет топологический изоморфизм
$$
(A,B):H^{s}(\Omega)\leftrightarrow
H^{s-2q}(\Omega)\times\prod_{j=1}^{q}H^{s-m_{j}-1/2}(\Gamma)\quad\mbox{при}
\;\;s\geq2q. \eqno(1.3)
$$
Здесь $H^{\sigma}(\Omega)$ и $H^{\sigma}(\Gamma)$, где
$\sigma\in\mathbb{R}$, --- гильбертовы пространства Соболева в
$\Omega$ и на $\Gamma$ соответственно.

Легко заметить, что этот результат не верен в случае произвольного
вещественного\;$s$. Так, при $s\leq m_{j}+1/2$ нельзя задать на
пространстве $H^{s}(\Omega)$ граничный дифференциальный оператор
$B_{j}$. Я. А. Ройтбергом \cite{Ro64}, \cite{Ro96}(п. 2.4) (см.
также \cite{Ber65}(гл. III, \S 6), \cite{Agr97}(п. 7.9))
предложено следующее определение обобщенного решения краевой
задачи (1.1), устраняющее этот недостаток.

В окрестности границы $\Gamma$ запишем дифференциальные выражения
$A$ и $B_{j}$ в виде
$$
A=\sum_{k=0}^{2q}\;A_{k}\,D_{\nu}^{k}\quad\mbox{и}\quad
B_{j}=\sum_{k=0}^{m_{j}}\;B_{j,k}\,D_{\nu}^{k}. \eqno(1.4)
$$
Здесь $D_{\nu}:=i\,\partial/\partial\nu$, где $\nu$ --- орт
внутренней нормали к границе $\Gamma$, а $A_{k}$ и $B_{j,k}$ ---
некоторые тангенциальные дифференциальные выражения. Интегрируя по
частям, запишем следующую формулу Грина:
$$
(Au,v)_{\Omega}=(u,A^{+}v)_{\Omega}-
i\sum_{k=1}^{2q}\;(D_{\nu}^{k-1}u,A^{(k)}v)_{\Gamma} \eqno(1.5)
$$
для произвольных функций $u,v\in
C^{\infty}(\,\overline{\Omega}\,)$. Здесь
$A^{(k)}:=\sum_{r=k}^{2q}D_{\nu}^{r-k}A_{r}^{+}$, где $A_{r}^{+}$
--- дифференциальное выражение, сопряженное к $A_{r}$. С помощью
предельного перехода получаем, что формула (1.5) справедлива для
каждого распределения $u\in H^{2q}(\Omega)$. Обозначим:
$$
u_{0}:=u\quad\mbox{и}\quad
u_{k}:=(D_{\nu}^{k-1}u)\upharpoonright\Gamma\quad\mbox{при}\quad
k=1,\ldots2q. \eqno(1.6)
$$

В силу (1.4), (1.5) краевая задача (1.1) относительного искомой
функции $u\in H^{2q}(\Omega)$ равносильна системе условий
$$
(u_{0},A^{+}v)_{\Omega}-
i\sum_{k=1}^{2q}\;(u_{k},A^{(k)}v)_{\Gamma}=(f,v)_{\Omega}
\quad\mbox{для любого}\quad v\in
C^{\infty}(\,\overline{\Omega}\,), \eqno(1.7)
$$
$$
\sum_{k=0}^{m_{j}}\;B_{j,k}\,u_{k+1}=g_{j}
\;\;\text{на}\;\;\Gamma\;\;\text{при}\;\;j=1,\ldots,q. \eqno(1.8)
$$
Заметим, что эти условия имеют смысл в случае произвольных (вообще
говоря, нерегулярных) распределений
$$
u_{0}\in\mathcal{D}'(\mathbb{R}^{n}),\;\;
\mathrm{supp}\,u_{0}\subseteq\overline{\Omega},\quad
u_{1},\ldots,u_{q}\in\mathcal{D}'(\Gamma). \eqno(1.9)
$$
Здесь, как обычно, через $\mathcal{D}'(\mathbb{R}^{n})$ и
$\mathcal{D}'(\Gamma)$ обозначены линейные топологические
пространства Шварца распределений в $\mathbb{R}^{n}$ и на $\Gamma$
соответственно. Поэтому введем следующее определение.

Вектор $u=(u_{0},u_{1},\ldots,u_{2q})$, удовлетворяющий условию
(1.9), называется \textit{обобщенным (по Ройтбергу) решением}
краевой задачи (1.1), если выполняются условия (1.7), (1.8).

Мы будем изучать обобщенные решения задачи  (1.1) в специально
подобранных парах гильбертовых пространств, построенных на основе
семейства пространств
$$
\left\{\,H^{s,\varphi}(\mathbb{R}^{n}):
\,s\in\mathbb{R},\varphi\in\mathcal{M}\right\}.
$$
Оно изучено авторами в \cite{MM2b} и названо \textit{уточненной}
шкалой в $\mathbb{R}^{n}$. Определение пространства
$H^{s,\varphi}(\mathbb{R}^{n})$ приведено ниже в п. 2. Здесь
отметим лишь, что это пространство гильбертово и состоит из
распределений в $\mathbb{R}^{n}$, гладкость которых
охарактеризована с помощью двух параметров --- числового $s$ и
функционального $\varphi$. Последний пробегает достаточно широкое
множество $\mathcal{M}$, состоящее из медленно меняющихся (по
Карамата) на $+\infty$ функций, и уточняет основную (степенную)
гладкость, задаваемую параметром $s$. В частном случае
$\varphi\equiv1$ пространство $H^{s,\varphi}(\mathbb{R}^{n})$
совпадает с пространством Соболева $H^{s}(\mathbb{R}^{n})$.

Пусть $s\in\mathbb{R}$, $\varphi\in\mathcal{M}$. В случае $s\geq0$
обозначим через $H^{s,\varphi,(0)}(\Omega)$ гильбертово
пространство сужений в область $\Omega$ всех распределений из
$H^{s,\varphi}(\mathbb{R}^{n})$. Далее, в случае $s<0$ обозначим
через $H^{s,\varphi,(0)}(\Omega)$ пространство, сопряженное к
пространству $H^{-s,\,1/\varphi,(0)}(\Omega)$ относительно
полуторалинейной формы $(\cdot,\cdot)_{\Omega}$. (Здесь уместно
отметить, что
$\varphi\in\mathcal{M}\Leftrightarrow1/\varphi\in\mathcal{M}$.)
Кроме того, обозначим через $H^{s,\varphi}(\Gamma)$ гильбертово
пространство распределений на $\Gamma$, принадлежащих локально
пространству $H^{s,\varphi}(\mathbb{R}^{n})$. (Детально указанные
пространства будут определены ниже в п. 2). Для каждого
$s\in\mathbb{R}\setminus\{1/2,\,3/2,\ldots,2q-1/2\}$ положим
$$
K_{s,\varphi,(2q)}(\Omega,\Gamma):=\bigl\{\,(u_{0},u_{1},\ldots,u_{2q}):\;
u_{0}\in H^{s,\varphi,(0)}(\Omega),\;u_{k}\in
H^{s-k+1/2,\varphi}(\Gamma),\;k=1,\ldots2q,
$$
$$
\mbox{причем}\;u_{k}=(D_{\nu}^{k-1}u)\upharpoonright\Gamma\;\;\mbox{если}\;\;
s>k-1/2\,\bigr\}. \eqno(1.10)
$$

Сформулируем основной результат статьи.

\textbf{Теорема 1.1.} \it В предположении $N=N^{+}=\{0\}$ оператор
$(A,B)$, соответствующий задаче \rm (1.1)\it, определяет
топологический изоморфизм
$$
(A,B):\,K_{s,\varphi,(2q)}(\Omega,\Gamma)\leftrightarrow
H^{s-2q,\varphi,(0)}(\Omega)\times\prod_{j=1}^{q}\,
H^{s-m_{j}-1/2,\varphi}(\Gamma) \eqno(1.11)
$$
для произвольных параметров
$s\in\mathbb{R}\setminus\{1/2,\,3/2,\ldots,2q-1/2\}$ и
$\varphi\in\mathcal{M}$. При этом решение
$u=(u_{0},u_{1},\ldots,u_{2q})$ задачи \rm (1.1) \it понимается
как обобщенное. \rm

Отождествляя функцию $u\in C^{\infty}(\,\overline{\Omega}\,)$ с
вектором $(u_{0},u_{1},\ldots,u_{2q})$, компоненты которого
вычисляются согласно (1.6), получаем, что операторы (1.3) и (1.11)
совпадают на множестве классических решений $u\in
C^{\infty}(\,\overline{\Omega}\,)$ краевой задачи (1.1). Это
множество плотно в пространствах $H^{s}(\Omega)$ и
$K_{s,\varphi,(2q)}(\Omega,\Gamma)$, являющихся областями
определения указаных операторов.

Более общее утверждение, чем теорема 1.1 приведено и доказано в п.
5.

\vspace{0.5cm}

\textbf{2. Уточненные шкалы пространств.} Сначала дадим
определение уточненной шкалы в $\mathbb{R}^{n}$, где
$n\in\mathbb{N}$ (см. \cite{MM2b}). Обозначим через $\mathcal{M}$
множество всех функций $\varphi:[1,+\infty)\rightarrow(0,+\infty)$
таких, что:

а) $\varphi$ измерима по Борелю на полуоси $[1,+\infty)$;

б) функции $\varphi$ и $1/\varphi$ ограничены на каждом отрезке
$[1,b]$, где $1<b<+\infty$;

в) функция $\varphi$ медленно меняющаяся на $+\infty$ по Карамата,
т. е. \cite{Se85}(п. 1.1)
$$
\lim_{t\rightarrow\,+\infty}\;\frac{\varphi(\lambda\,t)}{\varphi(t)}=1
\quad\mbox{для любого}\quad\lambda>0.
$$

Пусть $s\in\mathbb{R}$, $\varphi\in\mathcal{M}$. Обозначим через
$H^{s,\varphi}(\mathbb{R}^{n})$ пространство всех медленно
растущих распределений $u\in\mathcal{D}'(\mathbb{R}^{n})$ таких,
что преобразование Фурье $\widehat{u}$ распределения $u$ является
локально суммируемой по Лебегу в $\mathbb{R}^{n}$ функцией,
удовлетворяющей условию
$$
\int\langle\xi\rangle^{2s}\varphi^{2}(\langle\xi\rangle)\,|\widehat{u}(\xi)|^{2}
\,d\xi<\infty.
$$
Здесь интеграл берется по $\mathbb{R}^{n}$, а
$\langle\xi\rangle=(1+\xi_{1}^{2}+\ldots+\xi_{n}^{2})^{1/2}$ ---
сглаженный модуль вектора $\xi=(\xi_{1},\ldots,\xi_{n})\in
\mathbb{R}^{n}$. В пространстве $H^{s,\varphi}(\mathbb{R}^{n})$ в
качестве скалярного произведения возьмем величину
$$
\bigl(u,v\bigr)_{\mathrm{H}^{s,\varphi}(\mathbb{R}^{n})}:=
\int\langle\xi\rangle^{2s}\varphi^{2}(\langle\xi\rangle)
\,\widehat{u}(\xi)\,\overline{\widehat{v}(\xi)}d\xi.
$$
Она естественным образом порождает норму. Отметим, что мы
рассматриваем распределения, являющиеся \textit{антилинейными}
функционалами.

Пространство $H^{s,\varphi}(\mathbb{R}^{n})$ --- это частный
изотропный гильбертов случай пространств, рассмотренных Л.
Хермандером \cite[с. 54]{Her65}, \cite[с. 13]{Her86} и Б. П.
Волевичем, Л. Р. Панеяхом \cite[с. 14]{VoPa65}, \cite[с.
45]{Pa00}. В частном случае $\varphi\equiv1$ пространство
$H^{s,\varphi}(\mathbb{R}^{n})=H^{s,1}(\mathbb{R}^{n})$ совпадает
с пространством Соболева $H^{s}(\mathbb{R}^{n})$ порядка $s$. В
общем случае справедливы включения
$$
\bigcup_{\varepsilon>0}H^{s+\varepsilon}(\mathbb{R}^{n})=:H^{s+}(\mathbb{R}^{n})
\subset H^{s,\varphi}(\mathbb{R}^{n})\subset
H^{s-}(\mathbb{R}^{n}):=\bigcap_{\varepsilon>0}H^{s-\varepsilon}(\mathbb{R}^{n}).
\eqno(2.1)
$$
Они означают, что в семействе гильбертовых сепарабельных
пространств
$$
\{H^{s,\varphi}(\mathbb{R}^{n}):s\in\mathbb{R},\varphi\in\mathcal{M}\,\}
$$
функциональный параметр $\varphi$ \textit{уточняет} основную
(степенную) $s$-гладкость. Поэтому это семейство естественно
назвать \textit{уточненной} шкалой в $\mathbb{R}^{n}$ (по
отношению к соболевской шкале).

Теперь, следуя стандартной процедуре, определим аналоги
пространства $H^{s,\varphi}(\mathbb{R}^{n})$ для областей
$\overline{\Omega}$ и $\Omega$ (см. \cite{MM6}).

Обозначим
$$
H^{s,\varphi}_{\overline{\Omega}}(\mathbb{R}^{n}):=\left\{u\in
H^{s,\varphi}(\mathbb{R}^{n}):\,\mathrm{supp}\,u\subseteq
\overline{\Omega}\right\}.
$$
Отметим, что $H^{s,\varphi}_{\overline{\Omega}}(\mathbb{R}^{n})$
--- гильбертово сепарабельное пространство относительно скалярного
произведения в $H^{s,\varphi}(\mathbb{R}^{n})$.

Далее, положим
$$
H^{s,\varphi}(\Omega):=\bigl\{u\upharpoonright\Omega:\,u\in
H^{s,\varphi}(\mathbb{R}^{n})\bigr\},\quad
\bigl\|\,v\,\bigr\|_{H^{s,\varphi}(\Omega)}:=
\inf\,\bigl\{\,\bigl\|\,u\,\bigr\|_{H^{s,\varphi}(\mathbb{R}^{n})}:\,
u=v\;\;\mbox{в}\;\;\Omega\,\bigr\}.
$$
Пространство $H^{s,\varphi}(\Omega)$ сепарабельное и гильбертово,
поскольку норма в нем порождена скалярным произведением
$$
\bigl(v_{1},v_{2}\bigr)_{H^{s,\varphi}(\Omega)}:= \bigl(u_{1}-\Pi
u_{1},u_{2}-\Pi u_{2}\bigr)_{H^{s,\varphi}(\mathbb{R}^{n})}.
$$
Здесь $u_{j}\in H^{s,\varphi}(\mathbb{R}^{n})$, $u_{j}=v_{j}$ в
$\Omega$, где $j=1,\,2$, а $\Pi$ --- ортопроектор пространства
$H^{s,\varphi}(\mathbb{R}^{n})$ на подпространство $\{u\in
H^{s,\varphi}(\mathbb{R}^{n}):\mathrm{supp}\,u\subseteq\mathbb{R}^{n}
\setminus\Omega\}$.

Таким образом, пространство
$H^{s,\varphi}_{\overline{\Omega}}(\mathbb{R}^{n})$ состоит из
распределений, сосредоточенных в замкнутой области
$\overline{\Omega}$, а пространство $H^{s,\varphi}(\Omega)$ --- из
распределний, заданных в открытой области $\Omega$. Отметим
следующие их свойства \cite{MM6}(теорема 3.2). Множество
$$
C^{\infty}_{0}(\Omega):=\{u\in
C^{\infty}(\mathbb{R}^{n}):\mathrm{supp}\,u\subset\Omega\}
$$
плотно в $H^{s,\varphi}_{\overline{\Omega}}(\mathbb{R}^{n})$, а
множество $C^{\infty}(\,\overline{\Omega}\,)$ плотно в
$H^{s,\varphi}(\Omega)$. Пространства
$H^{s,\varphi}_{\overline{\Omega}}(\mathbb{R}^{n})$ и
$H^{-s,\,1/\varphi}(\Omega)$ взаимно сопряжены с равенством норм
относительно расширения по непрерывности полуторалинейной формы
$(u,v)_{\Omega}$, где $u\in C^{\infty}_{0}(\Omega)$, $v\in
C^{\infty}(\,\overline{\Omega}\,)$. Заметим здесь, что
пространство $H^{-s,\,1/\varphi}(\Omega)$ определено, поскольку
$1/\varphi\in\mathcal{M}$.

Рассмотрим также пространства распределений на многообразии
$\Gamma$. Возьмем конечный атлас из $C^{\infty}$-структуры на
$\Gamma$, образованный локальными картами
$\alpha_{j}:\mathbb{R}^{n}\leftrightarrow U_{j}$, где
$j=1,\ldots,k$. Здесь открытые множества $U_{j}$ составляют
конечное покрытие многообразия $\Gamma$. Пусть функции
$\chi_{j}\in C^{\infty}(\Gamma)$, где $j=1,\ldots,k$, образуют
разбиение единицы на $\Gamma$, удовлетворяющее условию
$\mathrm{supp}\,\chi_{j}\subset U_{j}$. Положим
$$
H^{s,\varphi}(\Gamma):=\bigl\{\,g\in\mathcal{D}'(\Gamma):\,
(\chi_{j}\,g)\circ\alpha_{j}\in
H^{s,\varphi}(\mathbb{R}^{n-1})\;\;\mbox{для каждого}\;\;
j=1,\ldots,k\bigr\},
$$
$$
\bigl(g_{1},g_{2}\bigr)_{H^{s,\varphi}(\Gamma)}:=
\sum_{j=1}^{k}\;\bigl(\,(\chi_{j}\,g_{1})\circ\alpha_{j},
(\chi_{j}\,g_{2})\circ\alpha_{j}\,\bigr)_{H^{s,\varphi}(\mathbb{R}^{n-1})}.
\eqno(2.2)
$$
Здесь $h\circ\alpha_{j}$ --- представление распределения
$h\in\mathcal{D}'(\Gamma)$ в локальной карте $\alpha_{j}$.
Скалярное произведение (2.2) естественным образом порождает норму
в пространстве $H^{s,\varphi}(\Gamma)$. Это пространство
гильбертово сепарабельное и с точностью до эквивалентных норм не
зависит от выбора атласа и разбиения единицы \cite{MM2b}(п. 3).
Множество $C^{\infty}(\Gamma)$ плотно в $H^{s,\varphi}(\Gamma)$.

Отметим далее следующее \cite{MM2b} (п. 3). Если $s>1/2$, то для
каждой функции $u\in H^{s,\varphi}(\Omega)$ определен по замыканию
ее след на границе $\Gamma$ --- функция $u\upharpoonright\Gamma\in
H^{s-1/2,\varphi}(\Gamma)$. При этом
$$
H^{s-1/2,\varphi}(\Gamma)=\{u\upharpoonright\Gamma:\,u\in
H^{s,\varphi}(\Omega)\},\quad
\bigl\|\,g\,\bigr\|_{H^{s-1/2,\varphi}(\Gamma)}\asymp
\inf\left\{\,\bigl\|\,u\,\bigr\|_{H^{s,\varphi}(\Omega)}:
\,u\upharpoonright\Gamma=g\right\}.
$$
В случае $s<1/2$ нельзя корректно определить след произвольного
распределения $u\in H^{s,\varphi}(\Omega)$ на границе $\Gamma$.
Вводимые ниже пространства $H^{s,\varphi,(r)}(\Omega)$, где
$r\in\mathbb{N}$, лишены этого недостатка.

Определим для каждого целого $r\geq0$ шкалу пространств
$$
\{H^{s,\varphi,(r)}(\Omega):s\in\mathbb{R},\varphi\in\mathcal{M}\,\}.
\eqno(2.3)
$$
Она сыграет центральную роль при изучении эллиптической краевой
задачи (1.1).

Пусть сначала $r=0$. В случае $s\geq0$ обозначим через
$H^{s,\varphi,(0)}(\Omega)$ гильбертово пространство
$H^{s,\varphi}(\Omega)$. В случае $s<0$ обозначим через
$H^{s,\varphi,(0)}(\Omega)$ гильбертово пространство
$H^{s,\varphi}_{\overline{\Omega}}(\mathbb{R}^{n})$, сопряженное к
пространству $H^{-s,\,1/\varphi}(\Omega)$ относительно
полуторалинейной формы  $(\cdot,\cdot)_{\Omega}$.

С точки зрения приложений к дифференциальным операторам удобна
трактовка пространства $H^{s,\varphi,(0)}(\Omega)$ как пополнения
линеала $C^{\infty}(\,\overline{\Omega}\,)$ по соответствующей
норме. Действительно, пространство $H^{s,\varphi,(0)}(\Omega)$ при
$s\geq0$ совпадает с пополнением линеала
$C^{\infty}(\,\overline{\Omega}\,)$ по норме пространства
$H^{s,\varphi}(\Omega)$. Далее заметим, что естественно
отождествлять функции из пространства
$L_{2}(\Omega)=H^{0}(\Omega)$ с их продолжениями нулем в
$\mathbb{R}^{n}$; в этом смысле
$L_{2}(\Omega)=H^{0}_{\overline{\Omega}}(\mathbb{R}^{n})$. При
таком отождествлении получим в силу (2.1) включения
$$
C^{\infty}_{0}(\Omega)\subset
C^{\infty}(\,\overline{\Omega}\,)\subset
L_{2}(\Omega)=H^{0}_{\overline{\Omega}}(\mathbb{R}^{n})\subset
H^{s,\varphi}_{\overline{\Omega}}(\mathbb{R}^{n})=H^{s,\varphi,(0)}(\Omega)
\;\;\mbox{при}\;\;s<0.
$$
Отсюда вытекает, что функции класса
$C^{\infty}(\,\overline{\Omega}\,)$ (продолженные нулем в
$\mathbb{R}^{n}$) образуют плотный линеал в пространстве
$H^{s,\varphi,(0)}(\Omega)$ при $s<0$. Значит, это пространство
является пополнением множества функций $u\in
C^{\infty}(\,\overline{\Omega}\,)$ по норме
$$
\sup\left\{\,\frac{\left|\,\bigl(u,v\bigr)_{\Omega}\,\right|}
{\;\quad\quad\quad\bigl\|\,v\,\bigr\|_{H^{-s,1/\varphi}(\Omega)}}\,:
\,v\in H^{-s,1/\varphi}(\Omega),\,v\neq0\right\}.
$$

Пусть теперь $r\in\mathbb{N}$. Положим
$E_{r}:=\{k-1/2:k=1,\ldots,r\}$. В случае $s\in\mathbb{R}\setminus
E_{r}$ обозначим через $H^{s,\varphi,(r)}(\Omega)$ пополнение
линейного пространства $C^{\infty}(\,\overline{\Omega}\,)$ по
норме
$$
\bigl\|\,u\,\bigr\|_{H^{s,\varphi,(r)}(\Omega)}:=
\left(\bigl\|\,u\,\bigr\|_{H^{s,\varphi,(0)}(\Omega)}^{2}+
\sum_{k=1}^{r}\;\left\|(D_{\nu}^{k-1}u)\upharpoonright\Gamma\right\|
_{H^{s-k+1/2,\varphi}(\Gamma)}^{2}\right)^{1/2}.
$$
Эта норма гильбертова, следовательно, и пространство
$H^{s,\varphi,(r)}(\Omega)$ гильбертово. Оно сепарабельно, как
будет показано ниже в п. 4.

В случае $s\in E_{r}$ определим пространство
$H^{s,\varphi,(r)}(\Omega)$ посредством интерполяции:
$$
H^{s,\varphi,(r)}(\Omega):=
\bigl[\,H^{s-\varepsilon,\varphi,(r)}(\Omega),
H^{s+\varepsilon,\varphi,(r)}(\Omega)\,\bigr]_{1/2}
\quad\mbox{при}\quad0<\varepsilon<1. \eqno(2.4)
$$
Определение использованного здесь интерполяционного метода
приведено в п. 3. Далее в п. 7 будет показано, что пространство
(2.4) с точностью до эквивалентных норм не зависит от
$\varepsilon$.

Семейство гильбертовых сепарабельных пространств (2.3) называем
\textit{модифицированной} по Ройтбергу уточненной шкалой порядка
$r$. В случае $\varphi\equiv1$ (пространства Соболева) эта шкала
введена и изучена Я. А. Ройтбергом  \cite{Ro64}, \cite{Ro96}(п.
2.4) (см. также \cite{Ber65}(гл. III, \S 6), \cite[с. 171]{Kr72},
\cite{Agr97}(п. 7.9)).

В силу определения модифицированной шкалы, оператор следа
$u\mapsto u\upharpoonright\Gamma$, $u\in
C^{\infty}(\,\overline{\Omega}\,)$, продолжается по непрерывности
до ограниченного оператора, действующего из пространства
$H^{s,\varphi,(r)}(\Omega)$ в пространство
$H^{s-1/2,\varphi}(\Gamma)$ при любых $s\in\mathbb{R}$,
$r\in\mathbb{N}$. Более того, для произвольного $u\in
H^{s,\varphi,(2q)}(\Omega)$ корректно определен по формулам (1.6)
посредством замыкания вектор $(u_{0},u_{1},\ldots,u_{2q})$.
Поэтому можно трактовать $u$ как обобщенное решение
$(u_{0},u_{1},\ldots,u_{2q})$ задачи (1.1).

\vspace{0.5cm}

\textbf{3. Интерполяция с функциональным параметром.} Интерполяция
с функциональным параметром пар гильбертовых пространств --- это
естественное обобщение классического интерполяционного метода
\cite[с. 21-23]{LM71}, \cite[c. 251]{Kr72} на случай, когда в
качестве параметра интерполяции вместо степенной берется более
общая функция. Приведем определение и некоторые свойства такой
интерполяции. Для наших целей достаточно ограничиться
сепарабельными гильбертовыми пространствами.

Упорядоченную пару $[X_{0},X_{1}]$ комплексных гильбертовых
пространств $X_{0}$ и $X_{1}$ будем называть \textit{допустимой},
если пространства $X_{0}$, $X_{1}$ сепарабельные и справедливо
непрерывное плотное вложение $X_{1}\hookrightarrow X_{0}$.

Пусть задана допустимая пара $X:=[X_{0},X_{1}]$ гильбертовых
пространств. Как известно \cite[c. 22]{LM71}, для $X$ существует
такой изометрический изоморфизм $J:X_{1}\leftrightarrow X_{\,0}$,
что $J$ является самосопряженным положительно определенным
оператором в пространстве $X_{\,0}$ с областью определения
$X_{1}$. Оператор $J$ называется \textit{порождающим} для пары
$X$, этот оператор определяется парой $X$ однозначно.

Обозначим через $\mathcal{B}$ множество всех функций, заданных
положительных и измеримых по Борелю на полуоси $(0,+\infty)$.
Пусть $\psi\in\mathcal{B}$. Поскольку спектр опрератора $J$
является подмножеством полуоси $(0,+\infty)$, в пространстве
$X_{0}$ определен как функция от $J$ оператор $\psi(J)$. Область
определения оператора $\psi(J)$ есть линейное многообразие,
плотное в $X_{0}$. Обозначим через $[X_{0},X_{1}]_{\psi}$ или,
короче, $X_{\psi}$ область определения оператора $\psi(J)$,
наделенную скалярным произведением графика
$$
(u,v)_{X_{\psi}}=(u,v)_{X_{0}}+(\psi(J)u,\psi(J)v)_{X_{0}}.
$$
Пространство $X_{\psi}$ гильбертово сепарабельное, причем
справедливо непрерывное плотное вложение $X_{\psi}\hookrightarrow
X_{0}$.

Функцию $\psi\in\mathcal{B}$ называем \textit{интерполяционным
параметром}, если для произвольных допустимых пар
$X=[X_{0},X_{1}]\,,Y=[Y_{0},Y_{1}]$ гильбертовых пространств и для
любого линейного отображения $T$, заданного на $X_{0}$,
выполняется следующее условие. Если при $j=0,\,1$ сужение
отображения $T$ на пространство $X_{j}$ является ограниченным
оператором $T:X_{j}\rightarrow Y_{j}$, то и сужение отображения
$T$ на пространство $X_{\psi}$ является ограниченным оператором
$T:X_{\psi}\rightarrow Y_{\psi}$.

Иными словами, функция $\psi$ является интерполяционным параметром
тогда и только тогда, когда отображение $X\mapsto X_{\psi}$
является интерполяционным функтором, заданным на категории
допустимых пар $X$ гильбертовых пространств (см. \cite{Tr80}, п.
1.2.2). В этом случае будем говорить, что \it пространство
$X_{\psi}$ получено в результате интерполяции пары $X$ с
функциональным параметром $\psi$. \rm

Классический результат \cite[c. 41]{LM71}, \cite[c. 250
--255]{Kr72} в теории интерполяции гильбертовых пространств
состоит в том, что степенная функция $\psi(t)=t^{\,\theta}$
порядка $\theta\in(0,1)$ является интерполяционным параметром. В
этом случае $\theta$ естественным образом выступает в качестве
числового параметра интерполяции и интерполяционное пространство
$X_{\psi}$ обозначается через $X_{\theta}$. Нам понадобится
следующий, более широкий чем степенной, класс интерполяционных
параметров \cite{MM2a} (теорема 2.1, лемма 2.1).

\textbf{Предложение 3.1.} \it Пусть функция $\psi\in\mathcal{B}$
ограничена на каждом отрезке $[a;b]$, где $0<a<b<+\infty$. Пусть,
кроме того, $\psi$
--- правильно меняющаяся на $+\infty$ по Карамата функция порядка
$\theta$, где $0<\theta<1$, т. е. \rm \cite{Se85}(п. 1.1) \it
$$
\lim_{t\rightarrow\,+\infty}\;\frac{\psi(\lambda\,t)}{\psi(t)}=
\lambda^{\theta}\quad\mbox{для любого}\quad \lambda>0.
$$
Тогда $\psi$ является интерполяционным параметром. При этом
справедливы непрерывные плотные вложения $X_{1}\hookrightarrow
X_{\psi}\hookrightarrow X_{0}$.\rm

Ниже будут использованы следующие свойства интерполяции.

\textbf{Предложение 3.2} (\cite{Sh74}, теорема 4). \it Пусть
задано конечное число допустимых пар $[X_{0}^{(k)},X_{1}^{(k)}]$
гильбертовых пространств, где $k=1,\ldots,m$. Тогда для любой
функции $\psi\in\mathcal{B}$ справедливо
$$
\left[\,\prod_{k=1}^{m}X_{0}^{(k)},\;\prod_{k=1}^{m}X_{1}^{(k)}\right]_{\psi}=
\prod_{k=1}^{m}\left[X_{0}^{(k)},\,X_{1}^{(k)}\right]_{\psi}
\quad\mbox{с равенством норм}.
$$\rm

\textbf{Предложение 3.3} (\cite{MM3}, теорема 2). \it Пусть
интерполяционные параметры $\zeta,\eta,\chi\in\mathcal{B}$
удовлетворяют следующему условию: для каждого числа
$\varepsilon>0$ существуют положительные числа
$c_{1}(\varepsilon)$, $c_{2}(\varepsilon)$ такие, что
$$
1\leq c_{1}(\varepsilon)\,\zeta(t)\leq c_{2}(\varepsilon)\,\eta(t)
\quad\mbox{и}\quad1\leq c_{1}(\varepsilon)\,\chi(t)
\quad\mbox{при}\quad t>\varepsilon.
$$
Тогда для произвольной допустимой пары $X$ гильбертовых
пространств справедливо равенство пространств
$[X_{\zeta},X_{\eta}]_{\chi}=X_{\psi}$ с эквивалентностью норм.
Здесь функция $\psi(t):=\zeta(t)\,\chi(\eta(t)/\zeta(t))$
аргумента $t>0$ является интерполяционным параметром.\rm

Напомним следующее определение. Линейный ограниченный оператор
$T:X\rightarrow Y$, где $X,Y$ --- банаховы пространства,
называется \textit{фредгольмовым}, если его ядро конечномерно, а
область значений $T(X)$ замкнута в $Y$ и имеет там конечную
коразмерность. \textit{Индексом} фредгольмового оператора $T$
называется число $\mathrm{ind}\,T=\dim\ker T-\dim(Y/\,T(X))$.

\textbf{Предложение 3.4} (\cite{Gey65}, предложение 5.2). \it
Пусть задано две допустимые пары $X=[X_{\,0},X_{1}]$ и
$Y=[Y_{\,0},Y_{1}]$ гильбертовых пространств. Пусть, кроме того,
на $X_{\,0}$ задано линейное отображение $T$, для которого
существуют ограниченные фредгольмовы операторы $T:X_{j}\rightarrow
Y_{j}$, где $j=0,\,1$, имеющие общее ядро $N$ и одинаковый индекс
$\kappa$. Тогда для произвольного интерполяционного параметра
$\psi\in\mathcal{B}$ ограниченный оператор $T:X_{\psi}\rightarrow
Y_{\psi}$ фредгольмов с ядром $N$, областью значений $Y_{\psi}\cap
T(X_{\,0})$ и тем же индексом $\kappa$. \rm

\vspace{0.5cm}

\textbf{4. Свойства модифицированной уточненной шкалы.} Сначала
изучим модифицированную шкалу (2.3) порядка $r=0$. Отметим
следующие ее свойства, установленные в \cite{MM6} (теорема 3.3).

\textbf{Предложение 4.1. } \it Пусть $s,\sigma\in\mathbb{R}$ и
$\varphi,\chi\in\mathcal{M}$. Тогда:

а) если $|s|<1/2$, то нормы в пространствах
$H^{s,\varphi}_{\overline{\Omega}}(\mathbb{R}^{n})$ и
$H^{s,\varphi}(\Omega)$ эквивалентны на плотном линейном
многообразии $C^{\infty}_{0}(\Omega)$, что означает следующее
равенство пространств с эквивалентностью норм в них:
$$
H^{s,\varphi,(0)}(\Omega)=
H^{s,\varphi}_{\overline{\Omega}}(\mathbb{R}^{n})=H^{s,\varphi}(\Omega)
\quad\mbox{при}\quad |s|<1/2; \eqno(4.1)
$$

б) пространства $H^{s,\varphi,(0)}(\Omega)$ и
$H^{-s,1/\varphi,(0)}(\Omega)$ взаимно сопряжены (при $s\neq0$ с
равенством норм, а при $s=0$ с эквивалентностью норм) относительно
полуторалинейной формы $(\cdot,\cdot)_{\Omega}$;

в) если $s<\sigma$, то справедливо компактное плотное вложение
$H^{\sigma,\chi,(0)}(\Omega)\hookrightarrow
H^{s,\varphi,(0)}(\Omega)$;

г) если $\varphi(t)\leq c\,\chi(t)$ при $t\gg1$ для некоторого
числа $c>0$, то справедливо непрерывное плотное вложение
$H^{s,\chi,(0)}(\Omega)\hookrightarrow H^{s,\varphi,(0)}(\Omega)$;
это вложение компактно, если $\varphi(t)/\chi(t)\rightarrow0$ при
$t\rightarrow+\infty$;

д) неравенство
$$
\int_{1}^{\,+\infty}\frac{d\,t}{t\,\varphi^{\,2}(t)}<\infty,
\eqno(4.2)
$$
равносильно непрерывности вложения
$H^{\rho+n/2,\varphi,(0)}(\Omega)\hookrightarrow
C^{\rho}(\,\overline{\Omega}\,)$, где целое $\rho\geq0$;
непрерывность такого вложения влечет его компактность. \rm

В связи с пунктами в), г) предложения 4.1 отметим следующее.
Плотное непрерывное вложение
$H^{\sigma,\chi,(0)}(\Omega)\hookrightarrow
H^{s,\varphi,(0)}(\Omega)$ понимается следующим образом.
Существует число $c>0$ такое, что
$$
\bigl\|\,u\,\bigr\|_{H^{s,\varphi,(0)}(\Omega)}\leq c\,
\bigl\|\,u\,\bigr\|_{H^{\sigma,\chi,(0)}(\Omega)}\quad\mbox{для
любого}\quad u\in C^{\infty}(\,\overline{\Omega}\,).
$$
Кроме того, тождественное отображение, заданное на плотном
линейном многообразии $C^{\infty}(\,\overline{\Omega}\,)$,
продолжается по непрерывности до ограниченного линейного
инъективного оператора $I:H^{\sigma,\chi,(0)}(\Omega)\rightarrow
H^{s,\varphi,(0)}(\Omega)$ (он называется оператором вложения).
Аналогично понимается плотное непрерывное вложение
$H^{\sigma,\chi,(r)}(\Omega)\hookrightarrow
H^{s,\varphi,(r)}(\Omega)$ (см. ниже).

Следующая теорема устанавливает тот факт, что каждое пространство
$H^{s,\varphi,(0)}(\Omega)$ может быть получено в результате
интерполяции пары соболевских пространств с подходящим
функциональным параметром.

\textbf{Теорема 4.1.} \it Пусть заданы функция
$\varphi\in\mathcal{M}$ и положительные числа
$\varepsilon,\delta$. Положим $\psi(t):=
t^{\varepsilon/(\varepsilon+\delta)}\varphi(t^{1/(\varepsilon+\delta)})$
при $t\geq1$ и $\psi(t):=\varphi(1)$ при $0<t<1$. Тогда:

а) функция $\psi\in\mathcal{B}$ является интерполяционным
параметром;

б) для каждого числа $s\in\mathbb{R}$ такого, что
$s-\varepsilon>-1/2$ или $s+\delta<1/2$, справедливо
$$
\left[H^{s-\varepsilon,1,(0)}(\Omega),
\,H^{s+\delta,1,(0)}(\Omega)\right]_{\psi}=
H^{s,\varphi,(0)}(\Omega)\quad\mbox{с эквивалентностью норм}.
$$\rm

\it\textbf{Доказательство.} Пункт а). \rm Непосредственно
проверяется, что функция $\psi\in\mathcal{B}$ удовлетворяет
условию предложения 3.1, где
$\theta=\varepsilon/(\varepsilon+\delta)\in(0,1)$. Следовательно,
она является интерполяционным параметром.

\textit{Пункт б).} Как установлено в \cite{MM2b} (теоремы 3.5,
3.7) и \cite{MM6} (теорема 3.1), для произвольного
$s\in\mathbb{R}$ справедливы следующие равенства пространств с
эквивалентностью норм в них:
$$
\left[H^{s-\varepsilon,1}(\Omega),H^{s+\delta,1}(\Omega)\right]_{\psi}
=H^{s,\varphi}(\Omega), \eqno (4.3)
$$
$$
\left[H^{s-\varepsilon,1}_{\overline{\Omega}}(\mathbb{R}^{n}),
H^{s+\delta,1}_{\overline{\Omega}}(\mathbb{R}^{n})\right]_{\psi}
=H^{s,\varphi}_{\overline{\Omega}}(\mathbb{R}^{n}). \eqno (4.4)
$$
Если $s-\varepsilon>-1/2$, то в силу (4.1) и (4.3) получаем
$$
\left[H^{s-\varepsilon,1,(0)}(\Omega),
\,H^{s+\delta,1,(0)}(\Omega)\right]_{\psi}=
\left[H^{s-\varepsilon,1}(\Omega),H^{s+\delta,1}(\Omega)\right]_{\psi}=
H^{s,\varphi}(\Omega)=H^{s,\varphi,(0)}(\Omega).
$$
Если  $s+\varepsilon<1/2$, то в силу (4.1) и (4.4) имеем
$$
\left[H^{s-\varepsilon,1,(0)}(\Omega),
\,H^{s+\delta,1,(0)}(\Omega)\right]_{\psi}=
\left[H^{s-\varepsilon,1}_{\overline{\Omega}}(\mathbb{R}^{n}),
H^{s+\delta,1}_{\overline{\Omega}}(\mathbb{R}^{n})\right]_{\psi}
=H^{s,\varphi}_{\overline{\Omega}}(\mathbb{R}^{n})=
H^{s,\varphi,(0)}(\Omega).
$$
Здесь наряду с равенством пространств выполняется эквивалентность
норм в них. Пункт б) доказан.

Далее изучим модифицированную шкалу (2.3) порядка
$r\in\mathbb{N}$. Нам понадобятся следующие свойства пространства
$H^{s,\varphi}(\Gamma)$, установленные в \cite{MM2b} (теоремы 3.5,
3.6, 3.8).

\textbf{Предложение 4.2. } \it Пусть $s,\sigma\in\mathbb{R}$ и
$\varphi,\chi\in\mathcal{M}$. Тогда:

а) для произвольных положительных чисел $\varepsilon,\delta$
справедливо
$$
\left[H^{s-\varepsilon,1}(\Gamma),
\,H^{s+\delta,1}(\Gamma)\right]_{\psi}=
H^{s,\varphi}(\Gamma)\quad\mbox{с эквивалентностью норм},
$$
где $\psi$ --- интерполяционный параметр из теоремы \rm 4.1\it;

б) пространства $H^{s,\varphi}(\Gamma)$ и
$H^{-s,1/\varphi}(\Gamma)$ взаимно сопряжены (с эквивалентностью
норм) относительно полуторалинейной формы
$(\cdot,\cdot)_{\Gamma}$;

в) пункты в), г) предложения \rm 4.1 \it сохраняют силу, если в их
формулировках заменить пространства $H^{\sigma,\chi,(0)}(\Omega)$,
$H^{s,\varphi,(0)}(\Omega)$, $H^{s,\chi,(0)}(\Omega)$ на
пространства $H^{\sigma,\chi}(\Gamma)$, $H^{s,\varphi}(\Gamma)$,
$H^{s,\chi}(\Gamma)$ соответственно;

г) неравенство \rm (4.2) \it равносильно непрерывности вложения
$H^{\rho+(n-1)/2,\varphi}(\Gamma)\hookrightarrow
C^{\rho}(\Gamma)$, где целое $\rho\geq0$; непрерывность такого
вложения влечет его компактность.

д) для любых $k\in\mathbb{N}$, $s>k-1/2$ линейное отображение
$u\mapsto(D_{\nu}^{k-1}u)\upharpoonright\Gamma$, где $u\in
C^{\infty}(\,\overline{\Omega}\,)$, продолжается по непрерывности
до ограниченного оператора, действующего из пространства
$H^{s,\varphi,(0)}(\Omega)=H^{s,\varphi}(\Omega)$ в пространство
$H^{s-k+1/2,\varphi}(\Gamma)$. \rm

Для произвольных $s\in\mathbb{R}$, $\varphi\in\mathcal{M}$ положим
$$
\Pi_{s,\varphi,(r)}(\Omega,\Gamma):=H^{s,\varphi,(0)}(\Omega)\times
\prod_{k=1}^{r}\,H^{s-k+1/2,\varphi}(\Gamma).
$$
Кроме того, если $s\notin E_{r}$, обозначим (см. (1.10))
$$
K_{s,\varphi,(r)}(\Omega,\Gamma):=
\bigl\{\,(u_{0},u_{1},\ldots,u_{r})\in\Pi_{s,\varphi,(r)}(\Omega,\Gamma):
$$
$$
u_{k}=(D_{\nu}^{k-1}u)\upharpoonright\Gamma\;\mbox{для
всех}\;k=1,\ldots r\;\mbox{таких, что}\;s>k-1/2\}.
$$
В силу предложения 4.2 д), $K_{s,\varphi,(r)}(\Omega,\Gamma)$
--- (замкнутое) подпространство в
$\Pi_{s,\varphi,(r)}(\Omega,\Gamma)$.

\textbf{Теорема 4.2. } \it Пусть $r\in\mathbb{N}$,
$s\in\mathbb{R}\setminus E_{r}$, $\varphi\in\mathcal{M}$. Тогда:

а) линейное отображение
$$
T_{r}:u\mapsto\left(\,u,u\upharpoonright\Gamma,\ldots,
(D_{\nu}^{r-1}u)\upharpoonright\Gamma\,\right),\quad u\in
C^{\infty}(\,\overline{\Omega}\,), \eqno(4.5)
$$
продолжается по непрерывности до изометрического изоморфизма
$$
T_{r}:\,H^{s,\varphi,(r)}(\Omega)\leftrightarrow
K_{s,\varphi,(r)}(\Omega,\Gamma); \eqno(4.6)
$$

б) для произвольных положительных чисел $\varepsilon,\delta$
таких, что числа $s,s-\varepsilon,s+\delta$ принадлежат одному из
интервалов
$$\alpha_{0}:=(-\infty,1/2),\;\;
\alpha_{k}:=(k-1/2,\,k+1/2),\;k=1,\ldots,r-1,\;\;
\alpha_{r}:=(r-1/2,+\infty),
$$
справедливо
$$
\left[H^{s-\varepsilon,1,(r)}(\Omega),
\,H^{s+\delta,1,(r)}(\Omega)\right]_{\psi}=
H^{s,\varphi,(r)}(\Omega)\quad\mbox{с эквивалентностью норм},
\eqno(4.7)
$$
где $\psi$ --- интерполяционный параметр из теоремы \rm 4.1\it.
\rm

\it\textbf{Доказательство.} \rm В случае $\varphi\equiv1$
(модификация соболевских пространств) пункт а) установлен Я. А.
Ройтбергом \cite{Ro96} (лемма 2.2.1). Выведем отсюда пункт б) для
произвольного $\varphi\in\mathcal{M}$, а затем пункт а).

Обозначим через $X_{\psi}$ левую часть равенства (4.7). Рассмотрим
изометрические операторы
$$
T_{r}:\,H^{\sigma,1,(r)}(\Omega)\rightarrow
\Pi_{\sigma,1,(r)}(\Omega,\Gamma),\quad\sigma\in\{s-\varepsilon,s+\delta\}.
$$
Применив к ним интерполяцию с параметром $\psi$, получим
ограниченный оператор
$$
T_{r}:\,X_{\psi}\rightarrow
\left[\,\Pi_{s-\varepsilon,1,(r)}(\Omega,\Gamma),\,
\Pi_{s+\delta,1,(r)}(\Omega,\Gamma)\,\right]_{\psi}. \eqno(4.8)
$$
В силу предложений 3.2, 4.2 а) и теоремы 4.1 имеем
$$
\left[\,\Pi_{s-\varepsilon,1,(r)}(\Omega,\Gamma),\,
\Pi_{s+\delta,1,(r)}(\Omega,\Gamma)\,\right]_{\psi}=
$$
$$
=\left[H^{s-\varepsilon,1,(0)}(\Omega),
\,H^{s+\delta,1,(0)}(\Omega)\right]_{\psi}\times
\prod_{k=1}^{r}\left[H^{s-\varepsilon-k+1/2,1}(\Gamma),
\,H^{s+\delta-k+1/2,1}(\Gamma)\right]_{\psi}=
$$
$$
=H^{s,\varphi,(0)}(\Omega)\times\prod_{k=1}^{r}\,H^{s-k+1/2,\varphi}(\Gamma)=
\Pi_{s,\varphi,(r)}(\Omega,\Gamma)
$$
с эквивалентностью норм. Следовательно,
$$
\bigl\|\,u\,\bigr\|_{H^{s,\varphi,(r)}(\Omega)}=
\bigl\|\,T_{r}\,u\,\bigr\|_{\Pi_{s,\varphi,(r)}(\Omega,\Gamma)}\leq
c_{1}\,\bigl\|\,u\,\bigr\|_{X_{\psi}}\quad\mbox{для всех}\quad
u\in C^{\infty}(\,\overline{\Omega}\,). \eqno(4.9)
$$
Здесь $c_{1}$ --- норма оператора (4.8).

Докажем неравенство, обратное к (4.9). По условию,
$s,s-\varepsilon,s+\delta\in\alpha_{p}$ для некоторого номера
$p\in\{0,1,\ldots,r\}$. Рассмотрим линейное отображение
$$
T_{r,p}:\,u\mapsto
\left(\,u,\,\{(D_{\nu}^{k-1}u)\upharpoonright\Gamma:\,p+1\leq
k\leq r\}\,\right),\quad u\in C^{\infty}(\,\overline{\Omega}\,).
$$
(Как и прежде, индекс $k$ целый.) Это отображение продолжается по
непрерывности до топологического изоморфизма
$$
T_{r,p}:\,H^{\sigma,1,(r)}(\Omega)\,\leftrightarrow\,
H^{\sigma,1,(0)}(\Omega)\,\times\prod_{p+1\leq k\leq
r}\,H^{\sigma-k+1/2,1}(\Gamma),\quad\sigma\in\{s-\varepsilon,s+\delta\}.
\eqno(4.10)
$$

Действительно, существование и ограниченность оператора (4.10)
следует из определения пространства $H^{\sigma,1,(r)}(\Omega)$.
Покажем, что этот оператор биективный. Пусть $u\in
H^{\sigma,1,(r)}(\Omega)$,
$$
\bigl(\,u_{0},\,\{u_{k}:\,p+1\leq k\leq r\}\,\bigr)\in
H^{\sigma,1,(0)}(\Omega)\,\times\prod_{p+1\leq k\leq
r}\,H^{\sigma-k+1/2,1}(\Gamma).
$$
Положим $u_{k}:=(D_{\nu}^{k-1}u_{0})\upharpoonright\Gamma$ при
$1\leq k\leq p$. Распределение $u_{k}$ определено корректно в силу
предложения 4.2 д), поскольку $\sigma>k-1/2$ для указанных номеров
$k$. Заметим, что $\sigma<k-1/2$ при $p+1\leq k\leq r$. Поэтому
$(u_{0},u_{1},\ldots,u_{r})\in K_{\sigma,1,(r)}(\Omega,\Gamma)$.
Как говорилось выше, пункт а) известен в случае $\varphi\equiv1$.
Значит, существуют топологические изоморфизмы
$$
T_{r}:H^{\sigma,1,(r)}(\Omega)\leftrightarrow
K_{\sigma,1,(r)}(\Omega,\Gamma),\quad\sigma\in\{s-\varepsilon,s+\delta\}.
$$
Отсюда, поскольку
$$
T_{r}\,u=(u_{0},u_{1},\ldots,u_{r})\;\Leftrightarrow\;
T_{r,p}\,u=(\,u_{0},\,\{u_{k}:\,p+1\leq k\leq r\}\,),
$$
вытекает, что ограниченный оператор (4.10) биективный.
Следовательно, по теореме Банаха об обратном операторе, (4.10) ---
топологический изоморфизм.

Применим к (4.10) интерполяцию с параметром $\psi$. В силу
предложений 3.2, 4.2 а) и теоремы 4.1 получим топологический
изоморфизм
$$
T_{r,p}:\,X_{\psi}\,\leftrightarrow\,
H^{s,\varphi,(0)}(\Omega)\,\times\prod_{p+1\leq k\leq
r}\,H^{s-k+1/2,\varphi}(\Gamma). \eqno(4.11)
$$
Отсюда вытекает неравенство, обратное к (4.9):
$$
\bigl\|\,u\,\bigr\|_{X_{\psi}}\leq c_{2}
\left(\bigl\|\,u\,\bigr\|_{H^{s,\varphi,(0)}(\Omega)}^{2}+
\sum_{p+1\leq k\leq
r}\;\left\|(D_{\nu}^{k-1}u)\upharpoonright\Gamma\right\|
_{H^{s-k+1/2,\varphi}(\Gamma)}^{2}\right)^{1/2}\leq c_{2}\,
\bigl\|\,u\,\bigr\|_{H^{s,\varphi,(r)}(\Omega)}
$$
для всех $u\in C^{\infty}(\,\overline{\Omega}\,)$. Здесь $c_{2}$
--- норма оператора, обратного к (4.11). Таким образом, нормы в
пространствах $X_{\psi}$ и $H^{s,\varphi,(r)}(\Omega)$
эквивалентны на множестве $C^{\infty}(\,\overline{\Omega}\,)$. Оно
плотно в $H^{s,\varphi,(r)}(\Omega)$ по определению и в $X_{\psi}$
в силу предложения 3.1. Следовательно,
$X_{\psi}=H^{s,\varphi,(r)}(\Omega)$ с точностью до эквивалентных
норм. Пункт б) доказан.

Докажем пункт а). Согласно определению пространства
$H^{s,\varphi,(r)}(\Omega)$ отображение (4.5) продолжается по
непрерывности до изометрического оператора
$$
T_{r}:\,H^{s,\varphi,(r)}(\Omega)\rightarrow
\Pi_{s,\varphi,(r)}(\Omega,\Gamma). \eqno(4.12)
$$
На основании предложения 4.2 д) справедливо включение
$T_{r}(H^{s,\varphi,(r)}(\Omega))\subseteq
K_{s,\varphi,(r)}(\Omega,\Gamma)$. Докажем обратное включение.
Пусть $(u_{0},u_{1},\ldots,u_{r})\in
K_{s,\varphi,(r)}(\Omega,\Gamma)$. В силу (4.11) и равенства
$X_{\psi}=H^{s,\varphi,(r)}(\Omega)$ существует топологический
изоморфизм
$$
T_{r,p}:\,H^{s,\varphi,(r)}(\Omega)\,\leftrightarrow\,
H^{s,\varphi,(0)}(\Omega)\,\times\prod_{p+1\leq k\leq
r}\,H^{s-k+1/2,\varphi}(\Gamma).
$$
Поэтому найдется такое $u\in H^{s,\varphi,(r)}(\Omega)$, что
$$
T_{r,p}\,u=(\,u_{0},\,\{u_{k}:\,p+1\leq k\leq r\}\,).
$$
Отсюда в силу предложения 4.2 д) вытекает равенство
$T_{r}\,u=(u_{0},u_{1},\ldots,u_{r})$. Тем самым, доказано
включение $K_{s,\varphi,(r)}(\Omega,\Gamma)\subseteq
T_{r}(H^{s,\varphi,(r)}(\Omega))$. Таким образом,
$T_{r}(H^{s,\varphi,(r)}(\Omega))=
K_{s,\varphi,(r)}(\Omega,\Gamma)$, что вместе с изометрическим
оператором (4.12) влечет изометрический изоморфизм (4.6). Пункт а)
доказан.

Теорема 4.2 доказана.

\textbf{Теорема 4.3.} \it Пусть $r\in\mathbb{N}$,
$s,\sigma\in\mathbb{R}$ и $\varphi,\chi\in\mathcal{M}$. Тогда:

а) гильбертово пространство $H^{s,\varphi,(r)}(\Omega)$
сепарабельное;

б) множество $C^{\infty}(\,\overline{\Omega}\,)$ плотно в
пространстве $H^{s,\varphi,(r)}(\Omega)$;

в) если $s>r-1/2$, то
$H^{s,\varphi,(r)}(\Omega)=H^{s,\varphi}(\Omega)$ с
эквивалентностью норм;

г) пункты в), г) предложения \rm 4.1 \it сохраняют силу, если в их
формулировках в обозначениях пространств заменить $(0)$ на $(r)$.
\rm

\it\textbf{Доказательство.} Пункт а). \rm Для $s\notin E_{r}$
сепарабельность пространства $H^{s,\varphi,(r)}(\Omega)$ вытекает
из теоремы 4.2 а) и сепарабельности пространства
$K_{s,\varphi,(r)}(\Omega,\Gamma)$. Если $s\in E_{r}$, то
пространство $H^{s,\varphi,(r)}(\Omega)$ сепарабельно в силу (2.4)
как результат интерполяции сепарабельных гильбертовых пространств.

\textit{Пункт б)} в случае $s\notin E_{r}$ содержится в
определеннии пространства $H^{s,\varphi,(r)}(\Omega)$. Если $s\in
E_{r}$, то в силу (2.4) и предложения 3.1 справедливо непрерывное
плотное вложение
$H^{s+\varepsilon,\varphi,(r)}(\Omega)\hookrightarrow
H^{s,\varphi,(r)}(\Omega)$ для достаточно малого $\varepsilon>0$.
Поскольку $s+\varepsilon\notin E_{r}$, множество
$C^{\infty}(\,\overline{\Omega}\,)$ плотно в пространстве
$H^{s+\varepsilon,\varphi,(r)}(\Omega)$. Следовательно, это
множество плотно и в пространстве $H^{s,\varphi,(r)}(\Omega)$.

\textit{Пункт в).} Если $s>r-1/2$, то в силу предложения 4.2 д)
нормы в пространствах $H^{s,\varphi,(r)}(\Omega)$ и
$H^{s,\varphi,(0)}(\Omega)=H^{s,\varphi}(\Omega)$ эквивалентны на
плотном линейном многообразии $C^{\infty}(\,\overline{\Omega}\,)$.
Следовательно, эти пространства равны.

\textit{Пункт г)} для $s,\sigma\notin E_{r}$ вытекает из
предложений 4.1 в), г) и 4.2 в) в силу теоремы 4.2 а). Если
$\{s,\sigma\}\cap E_{r}\neq\emptyset$ и $s<\sigma$, то в силу
(2.4) и предложения 3.1 для достаточно малого числа
$\varepsilon>0$ справедливы непрерывные плотные вложения
$$
H^{\sigma,\chi,(r)}(\Omega)\hookrightarrow
H^{\sigma-\varepsilon,\chi,(r)}(\Omega)\hookrightarrow
H^{s+\varepsilon,\varphi,(r)}(\Omega)\hookrightarrow
H^{s,\varphi,(r)}(\Omega).
$$
Наконец, если $s\in E_{r}$, то в силу (2.4) имеем
$$
H^{s\mp\varepsilon,\chi,(r)}(\Omega)\hookrightarrow
H^{s\mp\varepsilon,\varphi,(r)}(\Omega)\quad\Rightarrow\quad
H^{s,\chi,(r)}(\Omega)\hookrightarrow H^{s,\varphi,(r)}(\Omega),
$$
причем наследуется как непрерывность, так и компактность
\cite{Tr80} (п. 1.16.4) вложений.

Теорема 4.3 доказана.

\vspace{0.5cm}

\textbf{5. Эллиптическая краевая задача в модифицированной
уточненной шкале.} Напомним, что краевая задача (1.1) регулярная
эллиптическая, а $N$ и $N^{+}$ --- конечномерные бесконечно
гладкие ядра операторов задач (1.1) и (1.2) соответственно.

\textbf{Теорема 5.1.} \it Для произвольных параметров
$s\in\mathbb{R}$ и $\varphi\in\mathcal{M}$ линейное отображение
$$
(A,B):\,u\rightarrow(Au,B_{1}u,\ldots,B_{q}u),\quad u\in
C^{\infty}(\,\overline{\Omega}\,), \eqno(5.1)
$$
продолжается по непрерывности до ограниченного оператора
$$
(A,B):\,H^{s,\varphi,(2q)}(\Omega)\rightarrow
H^{s-2q,\varphi,(0)}(\Omega)\times\prod_{j=1}^{q}\,
H^{s-m_{j}-1/2,\varphi}(\Gamma)=:\mathcal{H}_{s,\varphi}(\Omega,\Gamma).
\eqno(5.2)
$$
Этот оператор фредгольмов. Его ядро совпадает с $N$, а область
значений равна множеству
$$
\left\{(f,g_{1},\ldots,g_{q})\in\mathcal{H}_{s,\varphi}(\Omega,\Gamma):\,
(f,v)_{\Omega}+\sum_{j=1}^{q}\,(g_{j},C^{+}_{j}v)_{\Gamma}=0\;\;
\mbox{для всех}\;\;v\in N^{+}\right\}. \eqno(5.3)
$$
Индекс оператора \rm (5.2) \it равен $\dim N-\dim N^{+}$ и не
зависит от $s,\varphi$. \rm

\it\textbf{Доказательство.} \rm В соболевском случае
$\varphi\equiv1$ эта теорема установлена Я. А. Ройтбергом
\cite{Ro96} (теоремы 4.1.1 и 5.3.1). Отсюда мы выведем общий
случай $\varphi\in\mathcal{M}$ с помощью интерполяции.

Сначала предположим, что $s\notin E_{2q}$. Пусть положительное
число $\varepsilon=\delta$ такое как в теореме 4.2 б). Отображение
(5.1) продолжается по непрерывности до ограниченных фредгольмовых
операторов
$$
(A,B):\,H^{s\mp\varepsilon,1,(2q)}(\Omega)\rightarrow
\mathcal{H}_{s\mp\varepsilon,1}(\Omega,\Gamma). \eqno(5.4)
$$
Они имеют общие ядро $N$, дефектное подпространство
$$
\left\{(v,C^{+}_{1}v,\ldots,C^{+}_{q}v):\,v\in N^{+}\right\}
\eqno(5.5)
$$
и индекс $\dim N-\dim N^{+}$. Применим к (5.4) интерполяцию с
функциональным параметром $\psi$ из теоремы 4.1. В силу
предложения 3.4 получим ограниченный фредгольмов оператор
$$
(A,B):\,\left[H^{s-\varepsilon,1,(2q)}(\Omega),
H^{s+\varepsilon,1,(2q)}(\Omega)\right]_{\psi}\rightarrow
\left[\mathcal{H}_{s-\varepsilon,1}(\Omega,\Gamma),
\mathcal{H}_{s+\varepsilon,1}(\Omega,\Gamma)\right]_{\psi}
$$
Он означает существование оператора (5.2), удовлетворяющего
условию настоящей теоремы. Это вытекает из теорем 4.1 б), 4.2 б) и
предложений 3.2, 4.2 а).

Предположим теперь, что $s\in E_{2q}$. Выберем произвольное число
$\varepsilon\in(0,1)$. Поскольку $s\mp\varepsilon\notin E_{2q}$,
существуют, как было доказано, ограниченные фредгольмовы операторы
$$
(A,B):\,H^{s\mp\varepsilon,\varphi,(2q)}(\Omega)\rightarrow
\mathcal{H}_{s\mp\varepsilon,\varphi}(\Omega,\Gamma),
$$
имеющие общие ядро $N$, дефектное подпространство (5.5) и индекс
$\dim N-\dim N^{+}$. Применив к этим операторам интерполяцию со
степенным параметром $t^{1/2}$, получим в силу предложения 3.4 и
формулы (2.4) ограниченный фредгольмов оператор
$$
(A,B):\,H^{s,\varphi,(2q)}(\Omega)\rightarrow
\left[\mathcal{H}_{s-\varepsilon,\varphi}(\Omega,\Gamma),
\mathcal{H}_{s+\varepsilon,\varphi}(\Omega,\Gamma)\right]_{1/2}.
\eqno(5.6)
$$
Он имеет те же ядро, дефектное подпространство и индекс.

Покажем, что
$$
\left[\mathcal{H}_{s-\varepsilon,\varphi}(\Omega,\Gamma),
\mathcal{H}_{s+\varepsilon,\varphi}(\Omega,\Gamma)\right]_{1/2}=
\mathcal{H}_{s,\varphi}(\Omega,\Gamma)\quad\mbox{с эквивалетностью
норм}. \eqno(5.7)
$$
Пусть число $\delta>0$ такое, что $s-\varepsilon-\delta>-1/2$ (это
возможно, поскольку $s\geq1/2$). На основании теоремы 4.1 и
предложений 3.2, 4.2 а) запишем
$$
\mathcal{H}_{s\mp\varepsilon,\varphi}(\Omega,\Gamma)=
\left[\mathcal{H}_{s-\varepsilon-\delta,1}(\Omega,\Gamma),
\mathcal{H}_{s+\varepsilon+\delta,1}(\Omega,\Gamma)\right]_{\psi_{\mp}}
\quad\mbox{с эквивалетностью норм}.
$$
Здесь интерполяционные параметры $\psi_{\mp}$ определяются по
формулам
$$
\psi_{-}(t):=t^{\delta/(2\varepsilon+2\delta)}
\varphi(t^{1/(2\varepsilon+2\delta)}),\quad
\psi_{+}(t):=t^{(2\varepsilon+\delta)/(2\varepsilon+2\delta)}
\varphi(t^{1/(2\varepsilon+2\delta)})\quad\mbox{при}\quad t\geq1
$$
и $\psi_{\mp}(t):=1$ при $0<t<1$. Отсюда в силу теоремы 3.3 о
повторной интерполяции получаем следующие равенства пространств с
эквивалентностью норм в них:
$$
\left[\mathcal{H}_{s-\varepsilon,\varphi}(\Omega,\Gamma),
\mathcal{H}_{s+\varepsilon,\varphi}(\Omega,\Gamma)\right]_{1/2}=
$$
$$
=\bigl[\,\left[\mathcal{H}_{s-\varepsilon-\delta,1}(\Omega,\Gamma),
\mathcal{H}_{s+\varepsilon+\delta,1}(\Omega,\Gamma)\right]_{\psi_{-}},\,
\left[\mathcal{H}_{s-\varepsilon-\delta,1}(\Omega,\Gamma),
\mathcal{H}_{s+\varepsilon+\delta,1}(\Omega,\Gamma)\right]_{\psi_{+}}\,\bigr]
_{1/2}=
$$
$$
=\left[\mathcal{H}_{s-\varepsilon-\delta,1}(\Omega,\Gamma),
\mathcal{H}_{s+\varepsilon+\delta,1}(\Omega,\Gamma)\right]_{\psi}.
\eqno(5.8)
$$
Здесь интерполяционный параметр $\psi$ определяется по формулам
$$
\psi(t):=\psi_{-}(t)\,(\psi_{+}(t)/\psi_{-}(t))^{1/2}=
t^{(\varepsilon+\delta)/(2\varepsilon+2\delta)}
\varphi(t^{1/(2\varepsilon+2\delta)})\quad\mbox{при}\quad t\geq1
$$
и $\psi(t)=1$ при $0<t<1$. Поэтому на основании тех же теоремы 4.1
и предложений 3.2, 4.2 а) имеем
$$
\left[\mathcal{H}_{s-\varepsilon-\delta,1}(\Omega,\Gamma),
\mathcal{H}_{s+\varepsilon+\delta,1}(\Omega,\Gamma)\right]_{\psi}=
\mathcal{H}_{s,\varphi}(\Omega,\Gamma)\quad\mbox{с эквивалетностью
норм}. \eqno(5.9)
$$
Теперь равенства (5.8), (5.9) влекут формулу (5.7).

В силу (5.7) ограниченный фредгольмов оператор (5.6) означает
существование оператора (5.2), удовлетворяющего условию настоящей
теоремы. Теорема 5.1 доказана.

Как отмечалось выше, теорема 5.1 уточняет применительно к шкале
пространств $H^{s,\varphi,(2q)}(\Omega)$ известный результат Я. А.
Ройтберга о свойствах эллиптической краевой задачи в
модифицированной шкале соболевских пространств \cite{Ro64},
\cite{Ro96} (см. также \cite[с. 169]{Kr72}, \cite[с. 248]{Ber65}).
В этой теореме $s$ --- \textit{произвольное} вещественное число.
Поэтому фредгольмовость оператора задачи установлена в
\textit{двусторонней} (иначе говоря, \textit{полной})
модифицированной уточненной шкале пространств. Заметим, что этот
оператор оставляет инвариантным параметр $\varphi\in\mathcal{M}$,
характеризюющий уточненную гладкость.

В силу теоремы 4.2 а) равенство $(A,B)u=f$, где $u\in
H^{s,\varphi,(2q)}(\Omega)$,
$f\in\mathcal{H}_{s,\varphi}(\Omega,\Gamma)$, равносильно тому,
что вектор $(u_{0},u_{1},\ldots,u_{2q}):=T_{2q}u$ является
обобщенным решением по Ройтбергу задачи (1.1). Указанный элемент
$u$ часто отождествляется с вектором $(u_{0},u_{1},\ldots,u_{2q})$
и также называется обобщенным решением задачи (1.1).

Из теоремы 4.3 в) вытекает, что оператор (5.2) совпадает с
оператором
$$
(A,B):\,H^{s,\varphi}(\Omega)\rightarrow
\mathcal{H}_{s,\varphi}(\Omega,\Gamma)\quad\mbox{при}\quad
s>2q-1/2.
$$
Его фредгольмовость установлена в \cite{MM2b} (теорема 4.1).

В частном случае $N=N^{+}=\{0\}$ оператор (5.2) является
топологическим изоморфизмом в силу теоремы 5.1 и теоремы Банаха об
обратном операторе. Следовательно, теорема 5.1 содержит теорему
1.1. В общем случае изоморфизм удобно задавать с помощью следующих
проекторов (ср. \cite{Ro96}, леммы 4.1.2 и 5.3.2).

\textbf{Лемма 5.1.} \it Для произвольных $s\in\mathbb{R}$ и
$\varphi\in\mathcal{M}$ существуют следующие разложения
пространств $H^{s,\varphi,(2q)}(\Omega)$ и
$\mathcal{H}_{s,\varphi}(\Omega,\Gamma)$ в прямые суммы замкнутых
подпространств:
$$
H^{s,\varphi,(2q)}(\Omega)=N\dotplus\left\{u\in
H^{s,\varphi,(2q)}(\Omega):\;(u_{0},w)_{\Omega}=0\;\;\mbox{для
любого}\;\;w\in N\right\}, \eqno(5.10)
$$
$$
\mathcal{H}_{s,\varphi}(\Omega,\Gamma)=
\left\{\,(v,0,\ldots,0):\;v\in N^{+}\right\}\dotplus
(A,B)\bigl(H^{s,\varphi,(2q)}(\Omega)\bigr). \eqno(5.11)
$$
Здесь $u_{0}$ --- начальная компонента вектора
$(u_{0},u_{1},\ldots,u_{2q}):=T_{2q}u$. Обозначим через $P$ косой
проектор пространства $H^{s,\varphi,(2q)}(\Omega)$ на второе
слагаемое суммы \rm (5.10)\it, а через $Q^{+}$ косой проектор
пространства $\mathcal{H}_{s,\varphi}(\Omega,\Gamma)$ на второе
слагаемое суммы \rm (5.11) \it (паралельно первому слагаемому).
Эти проекторы не зависят от $s,\varphi$. \rm

\it\textbf{Доказательство.} \rm Докажем (5.10). Из определения
пространства $H^{s,\varphi,(2q)}(\Omega)$ вытекает, что
отображение $u\mapsto u_{0}$ является ограниченным оператором
$T_{0}:H^{s,\varphi,(2q)}(\Omega)\rightarrow
H^{s,\varphi,(0)}(\Omega)$. Поэтому, второе слагаемое суммы (5.10)
--- замкнутое подпространство. Оно имеет тривиальное пересечение с $N$.
В силу предложения 4.1 б) и конечномерности подпространства $N$
справедливо разложение
$$
H^{s,\varphi,(0)}(\Omega)=N\dotplus\left\{u_{0}\in
H^{s,\varphi,(0)}(\Omega):\;(u_{0},w)_{\Omega}=0\;\;\mbox{для
любого}\;\;w\in N\right\}.
$$
Обозначим через $\Pi$ косой проектор на первое слагаемое этой
суммы параллельно второму слагаемому. Для произвольного $u\in
H^{s,\varphi,(2q)}(\Omega)$ запишем $u=u'+u''$, где $u':=\Pi
u_{0}\in N$, а $u'':=u-\Pi u_{0}\in H^{s,\varphi,(2q)}(\Omega)$
удовлетворяет условию $(u''_{0},w)_{\Omega}=(u_{0}-\Pi
u_{0},w)_{\Omega}=0$ при любом $w\in N$. Равенство (5.10)
доказано.

Равенство (5.11) вытекает из того, что в силу теоремы 5.1
подпространства, записанные в сумме (5.11), замкнутые, имеют
тривиальное пересечение и конечная размерность первого
пространства совпадает с коразмерностью второго. Независимость
проекторов $P$ и $Q^{+}$ от параметров $s,\varphi$ вытекает из
включений $N,N^{+}\subset C^{\infty}(\,\overline{\Omega}\,)$.
Лемма 5.1 доказана.

\textbf{Теорема 5.2.} \it Для произвольных параметров
$s\in\mathbb{R}$, $\varphi\in\mathcal{M}$ сужение оператора \rm
(5.2) \it на подпространство $P(H^{s,\varphi,(2q)}(\Omega))$
является топологическим изоморфизмом \rm
$$
(A,B):\,P\bigl(H^{s,\varphi,(2q)}(\Omega)\bigr)\,\leftrightarrow\,
Q^{+}\bigl(\mathcal{H}_{s,\varphi,(2q)}(\Omega,\Gamma)\bigr).
\eqno(5.12)
$$

\it\textbf{Доказательство.} \rm Согласно теореме 5.1, $N$ ---
ядро, а $Q^{+}(\mathcal{H}_{s,\varphi,(2q)}(\Omega,\Gamma))$
--- область значений оператора (5.2). Следовательно, оператор
(5.12) --- биекция. Кроме того, он ограничен. Значит, оператор
(5.12) является топологическим изоморфизмом в силу теоремы Банаха
об обратном операторе.

\it\textbf{Замечание} \rm \textbf{5.1.} Теорема 5.2 остается
верной, если заменить проектор $Q^{+}$ на оператор проектирования
$Q^{+}_{1}$ пространства
$\mathcal{H}_{s,\varphi,(2q)}(\Omega,\Gamma)$ на подпространство
$(A,B)(H^{s,\varphi,(2q)}(\Omega))$ паралельно дефектному
подпространству (5.5).

Из теоремы 5.2 вытекает следующая априорная оценка решения
эллиптической краевой задачи (1.1).

\textbf{Теорема 5.3.} \it Для произвольных фиксированных
параметров $s\in\mathbb{R}$, $\varphi\in\mathcal{M}$ и $\sigma<s$
существует число $c>0$ такое, что для каждого $u\in
H^{s,\varphi,(2q)}(\Omega)$ справедливо неравенство \rm
$$
\bigl\|\,u\,\bigr\|_{H^{s,\varphi,(2q)}(\Omega)}\leq c\,
\left(\,\bigl\|(A,B)\,u\,\bigr\|_{\mathcal{H}_{s,\varphi}(\Omega,\Gamma)}+
\bigl\|\,u\,\bigr\|_{H^{\sigma,1,(2q)}(\Omega)}\,\right).
\eqno(5.13)
$$

\it\textbf{Доказательство. }\rm Пусть $u\in
H^{s,\varphi,(2q)}(\Omega)$. Так как $N$ --- конечномерное
подпространство в пространствах $H^{s,\varphi,(2q)}(\Omega)$ и
$H^{\sigma,1,(2q)}(\Omega)$, нормы в этих пространствах
эквивалентны на $N$. В частности, для функции $u-Pu\in N$
справедливо неравенство
$$
\bigl\|\,u-Pu\,\bigr\|_{H^{s,\varphi,(2q)}(\Omega)}\leq\,
c_{1}\,\bigl\|\,u-Pu\,\bigr\|_{H^{\sigma,1,(2q)}(\Omega)}
$$
с постоянной $c_{1}>0$, не зависящей от $u$. Отсюда получаем
$$
\bigl\|\,u\,\bigr\|_{H^{s,\varphi,(2q)}(\Omega)}\leq\,
c_{1}\,\bigl\|\,u-Pu\,\bigr\|_{H^{\sigma,1,(2q)}(\Omega)}+
\bigl\|\,Pu\,\bigr\|_{H^{s,\varphi,(2q)}(\Omega)}\leq
$$
$$
\leq
c_{1}\,c_{2}\,\bigl\|\,u\,\bigr\|_{H^{\sigma,1,(2q)}(\Omega)}+
\bigl\|\,Pu\,\bigr\|_{H^{s,\varphi,(2q)}(\Omega)}, \eqno(5.14)
$$
где $c_{2}$ --- норма проектора $1-P$, действующего в пространстве
$H^{\sigma,1,(2q)}(\Omega)$. Далее, поскольку $(A,B)Pu=(A,B)u$, то
$Pu\in H^{s,\varphi,(2q)}(\Omega)$ --- прообраз распределения
$(A,B)u\in\mathcal{H}_{s,\varphi}(\Omega,\Gamma)$ при
топологическом изоморфизме (5.12). Следовательно,
$$
\bigl\|\,Pu\,\bigr\|_{H^{s,\varphi,(2q)}(\Omega)}\leq\,
c_{3}\,\bigl\|\,(A,B)u\,\bigr\|_{\mathcal{H}_{s,\varphi}(\Omega,\Gamma)},
$$
где $c_{3}$ --- норма оператора, обратного к (5.12). Отсюда и из
неравенства (5.14) вытекает оценка (5.13). Теорема 5.3 доказана.

\vspace{0.5cm}

\textbf{6. Локальная гладкость решения.} Предположим, что правые
части эллиптической краевой задачи (1.1) имеют на некотором
открытом в $\overline{\Omega}$ множестве дополнительную гладкость
в уточненной шкале пространств. Покажем, что обобщенное решение
$u$ унаследует такую же дополнительную гладкость на этом
множестве. Предварительно рассмотрим случай дополнительной
гладкости во всей области $\overline{\Omega}$.

\textbf{Теорема 6.1. }\it Пусть $s\in\mathbb{R}$. Предположим, что
элемент $u\in H^{s,1,(2q)}(\Omega)$ является обобщенным решением
задачи \rm (1.1)\it, где
$$
f\in H^{s+\varepsilon-2q,\varphi,(0)}(\Omega)\quad\mbox{и}\quad
g_{j}\in
H^{s+\varepsilon-m_{j}-1/2,\varphi}(\Gamma)\quad\mbox{при}\quad
j=1,\ldots,q
$$
для некоторых $\varepsilon\geq0$ и $\varphi\in\mathcal{M}$. Тогда
$u\in H^{s+\varepsilon,\varphi,(2q)}(\Omega)$.\rm

\it\textbf{Доказательство.} \rm  По условию и теореме 5.1 имеем
$$
F:=(f,g_{1},\ldots,g_{q})=(A,B)u\in(A,B)\bigl(H^{s,1,(2q)}(\Omega)\bigr)\cap
\mathcal{H}_{s+\varepsilon,\varphi}(\Omega,\Gamma)=
(A,B)\bigl(H^{s+\varepsilon,\varphi,(2q)}(\Omega)\bigr).
$$
Следовательно, существует такое $v\in
H^{s+\varepsilon,\varphi,(2q)}(\Omega)$, что $(A,B)v=F$. Отсюда
получаем $(A,B)(u-v)=0$, что в силу теоремы 5.1 влечет $w:=u-v\in
N\subset C^{\infty}(\,\overline{\Omega}\,)$. Таким образом,
поскольку $C^{\infty}(\,\overline{\Omega}\,)\subset
H^{s+\varepsilon,\varphi,(2q)}(\Omega)$, справедливо $u=v+w\in
H^{s+\varepsilon,\,\varphi}(\Omega)$, что и требовалось доказать.

Рассмотрим теперь случай локальной гладкости. Пусть $U$ ---
открытое непустое подмножество замкнутой области
$\overline{\Omega}$. Положим $\Omega_{0}:=U\cap\Omega$ и
$\Gamma_{0}:=U\cap\Gamma$ (возможен случай
$\Gamma_{0}=\emptyset$). Введем следующие локальные аналоги
пространств $H^{\sigma,\varphi,(r)}(\Omega)$ и
$H^{\sigma,\varphi}(\Gamma)$, где $\sigma\in\mathbb{R}$,
$\varphi\in\mathcal{M}$ и целое $r\geq0$. Положим
$$
H^{\sigma,\varphi,(r)}_{\mathrm{loc}}(\Omega_{0},\Gamma_{0}):=
$$
$$
:=\left\{u\in\bigcup_{s\in\mathbb{R}} H^{s,1,(r)}(\Omega):\,\chi
u\in H^{\sigma,\varphi,(r)}(\Omega)\;\mbox{для всех}\;\chi\in
C^{\infty}(\overline{\Omega}),\,\mathrm{supp}\,\chi\subset\Omega_{0}\cup\Gamma_{0}
\right\},
$$
$$
H^{\sigma,\varphi}_{\mathrm{loc}}(\Gamma_{0}):=\bigl\{h\in\mathcal{D}'(\Gamma):
\chi\,h\in H^{\sigma,\varphi}(\Gamma)\;\mbox{для всех}\;\chi\in
C^{\infty}(\Gamma),\,\mathrm{supp}\,\chi\subset \Gamma_{0}\bigr\}.
$$
В связи с определением пространства
$H^{\sigma,\varphi,(r)}_{\mathrm{loc}}(\Omega_{0},\Gamma_{0})$
отметим, что для произвольной функции $\chi\in
C^{\infty}(\,\overline{\Omega}\,)$ отображение $u\mapsto \chi u$,
$u\in C^{\infty}(\,\overline{\Omega}\,)$, продолжается по
непрерывности до ограниченного оператора в каждом пространстве
$H^{s,1,(r)}(\Omega)$ \cite{Ro96} (п. 2.3). Тем самым для $u\in
H^{s,1,(r)}(\Omega)$ корректно определено произведение $\chi u\in
H^{s,1,(r)}(\Omega)$.

\textbf{Теорема 6.2. }\it Пусть $s\in\mathbb{R}$. Предположим, что
элемент $u\in H^{s,1,(2q)}(\Omega)$ является обобщенным решением
задачи \rm (1.1)\it,  где
$$
f\in
H^{s+\varepsilon-2q,\varphi,(0)}_{\mathrm{loc}}(\Omega_{0},\Gamma_{0})
\quad\mbox{и}\quad g_{j}\in
H^{s+\varepsilon-m_{j}-1/2,\varphi}_{\mathrm{loc}}(\Gamma_{0})\quad\mbox{при}
\quad j=1,\ldots,q \eqno (6.1)
$$
для некоторых $\varepsilon\geq0$ и $\varphi\in\mathcal{M}$. Тогда
$u\in
H^{s+\varepsilon,\varphi,(2q)}_{\mathrm{loc}}(\Omega_{0},\Gamma_{0})$.\rm

\it\textbf{Доказательство.} \rm  Покажем сначала, что из условия
(6.1) вытекает следующее свойство повышения локальной гладкости
решения $u$: для каждого числа $r\geq1$ справедлива импликация
$$
u\in
H^{s+\varepsilon-r,\varphi,(2q)}_{\mathrm{loc}}(\Omega_{0},\Gamma_{0})
\;\Rightarrow\;u\in
H^{s+\varepsilon-r+1,\varphi,(2q)}_{\mathrm{loc}}(\Omega_{0},\Gamma_{0}).
\eqno(6.2)
$$

Выберем произвольно функции $\chi,\eta$ такие, что
$$
\chi,\eta\in C^{\infty}(\,\overline{\Omega}\,);\quad
\mathrm{supp}\,\chi,\,\mathrm{supp}\,\eta\subset\Omega_{0}\cup\Gamma_{0}
\quad\mbox{и}\quad\eta=1 \quad\mbox{в окрестности
}\quad\mathrm{supp}\,\chi. \eqno(6.3)
$$
Переставив оператор умножения на функцию $\chi$ с
дифференциальными операторами $A$ и $B_{j}$, где $j=1,\ldots,q$,
можно записать для произвольного $v\in
C^{\infty}(\,\overline{\Omega}\,)$ следующие равенства:
$$
(A,B)(\chi v)=(A,B)(\chi\eta v)=\chi(A,B)(\eta v)+(A',B')(\eta v)=
\chi(A,B)v+(A',B')(\eta v). \eqno(6.4)
$$
Здесь $A'$ --- некоторое линейное дифференциальное выражение в
$\overline{\Omega}$, а $B'=(B_{1}',\ldots,B_{q}')$ --- набор
граничных линейных дифференциальное выражений на $\Gamma$.
Коэффициенты этих выражений бесконечно гладкие, а порядки
удовлетворяют условиям $\mathrm{ord}\,A'\leq2k-1$ и
$\mathrm{ord}\,B'_{j}\leq m_{j}-1$. Отсюда следует, что
отображение $v\mapsto(A',B')v$, где $v\in
C^{\infty}(\,\overline{\Omega}\,)$, продолжается по непрерывности
до ограниченного оператора
$$
(A',B'):\,H^{\sigma,\varphi,(2q)}(\Omega)\rightarrow
\mathcal{H}_{\sigma+1,\varphi}(\Omega,\Gamma)\quad\mbox{для
произвольного}\quad\sigma\in\mathbb{R}. \eqno (6.5)
$$
В случае $\varphi\equiv1$ это доказано в \cite{Ro96} (п. 2.3).
Отсюда общий случай $\varphi\in\mathcal{M}$ выводится с помощью
интерполяции также как и в доказательстве теоремы 5.1. Аналогично
доказывается, что оператор умножения на функцию класса
$C^{\infty}(\,\overline{\Omega}\,)$ ограничен в пространствах
$H^{\sigma,\varphi,(r)}(\Omega)$ и
$\mathcal{H}_{\sigma,\varphi}(\Omega,\Gamma)$ для каждого
$\sigma\in\mathbb{R}$. Следовательно, равенство (6.4) продолжается
по непрерывности на класс функций $v\in
H^{\sigma,\varphi,(r)}(\Omega)$. Возьмем в этом равенстве $v:=u$,
где $u$ --- указанное в условии решение задачи (1.1). Запишем
$$
(A,B)(\chi u)=\chi F+(A',B')(\eta u). \eqno(6.6)
$$
Здесь вектор $F:=(f,g_{1},\ldots,g_{q})$ удовлетворяет в силу
(6.1) и (6.3) условию
$$
\chi F\in\mathcal{H}_{s+\varepsilon,\varphi}(\Omega,\Gamma).
\eqno(6.7)
$$

Предположим, что $u\in
H^{s+\varepsilon-r,\varphi,(2q)}_{\mathrm{loc}}(\Omega_{0},\Gamma_{0})$
для некоторого числа $r\geq1$. Тогда $\eta u\in
H^{s+\varepsilon-r,\varphi,(2q)}(\Omega)$, что вместе с формулами
(6.5) --  (6.7) влечет включение
$$
(A,B)(\chi u)
\in\mathcal{H}_{s+\varepsilon-r+1,\varphi}(\Omega,\Gamma).
$$
Отсюда в силу теоремы 6.1 следует свойство $\chi u\in
H^{s+\varepsilon-r+1,\varphi,(2q)}(\Omega)$. Оно ввиду
произвольности выбора функции $\chi$, удовлетворяющей условию
(6.3), означает включение $u\in
H^{s+\varepsilon-r+1,\varphi,(2q)}_{\mathrm{loc}}(\Omega_{0},\Gamma_{0})$.
Импликация (6.2) доказана.

Теперь легко вывести теорему из (6.2). В силу теоремы 4.3 г) имеем
$$
u\in H^{s,1,(2q)}(\Omega)\subset
H^{s+\varepsilon-k,\varphi,(2q)}(\Omega)\subseteq
H^{s+\varepsilon-k,\varphi,(2q)}_{\mathrm{loc}}(\Omega_{0},\Gamma_{0})
$$
для целого $k>\varepsilon$. Применив импликацию (6.2)
последовательно для значений $r=k,\ldots,1$, получим
$$
u\in
H^{s+\varepsilon-k,\varphi,(2q)}_{\mathrm{loc}}(\Omega_{0},\Gamma_{0})
\Rightarrow u\in
H^{s+\varepsilon-k+1,\varphi,(2q)}_{\mathrm{loc}}(\Omega_{0},\Gamma_{0})
\Rightarrow\ldots\Rightarrow u\in
H^{s+\varepsilon,\varphi,(2q)}_{\mathrm{loc}}(\Omega_{0},\Gamma_{0}),
$$
что и требовалось доказать.

В соболевском случае $\varphi\equiv1$ теоремы 6.1, 6.2 доказаны Я.
А. Ройтбергом \cite{Ro64}, \cite{Ro96} (гл. 7) (см. также
\cite{Ber65}, гл. III, \S 6).

В качестве приложения теорем 6.1 и 6.2 приведем одно достаточное
условие того, что обобщенное по Ройтбергу решение $u$
эллиптической краевой задачи (1.1) является классическим, т.е.
удовлетворяет условию
$$
u\in H^{\sigma+2q,1}(\Omega)\cap C^{2q}(\Omega)\cap
C^{m}(\,\overline{\Omega}\,), \eqno(6.8)
$$
где $\sigma>-1/2$, $m:=\max\{m_{1},\ldots,m_{q}\}$. Это условие
возникает следующим образом. В силу теоремы 4.3 в) и предложения
4.2 д) из включения
$$
u\in H^{\sigma+2q,1,(2q)}(\Omega)=H^{\sigma+2q,1}(\Omega)
$$
вытекает, что элемент $u$ является решением задачи (1.1) в смысле
теории распределений, заданных в области $\Omega$. Теперь
корректно рассматривать включение $u\in C^{2q}(\Omega)\cap
C^{m}(\,\overline{\Omega}\,)$. Оно означает, что в (1.1) функции
$Au$ и $B_{j}u$ вычисляются с помощью классических производных, т.
е. решение $u$ классическое.

\textbf{Теорема 6.3.} \it Пусть $s\in\mathbb{R}$. Предположим, что
элемент $u\in H^{s,1,(2q)}(\Omega)$ является обобщенным решением
задачи \rm (1.1)\it, где
$$
f\in H^{n/2,\varphi,(0)}_{\mathrm{loc}}(\Omega,\emptyset)\cap
H^{m-2q+n/2,\varphi,(0)}(\Omega)\cap H^{\sigma,1,(0)}(\Omega),
\eqno(6.9)
$$
$$
g_{j}\in H^{\,m-m_{j}+(n-1)/2,\,\varphi}(\Gamma)\cap
H^{\sigma+2q-m_{j}-1/2,\,1}(\Gamma)\;\;\mbox{при}\;\; j=1,\ldots,q
\eqno (6.10)
$$
для некоторых числа $\sigma>-1/2$ и функционального параметра
$\varphi\in\mathcal{M}$, удовлетворяющего неравенству \rm
(4.2)\it. Тогда решение $u$ классическое. \rm

\it\textbf{Доказательство.} \rm В силу теорем 6.1, 6.2 из условий
(6.9), (6.10) вытекает включение
$$
u\in H^{2q+n/2,\varphi,(2q)}_{\mathrm{loc}}(\Omega,\emptyset)\cap
H^{m+n/2,\varphi,(2q)}(\Omega)\cap H^{\sigma+2q,1,(2q)}(\Omega).
$$
Отсюда на основании предложения 4.1 д) и теоремы 4.3 в) имеем:
$$
u\in H^{m+n/2,\varphi,(2q)}(\Omega)\cap
H^{\sigma+2q,1,(2q)}(\Omega)= H^{m+n/2,\varphi}(\Omega)\cap
H^{\sigma+2q,1}(\Omega)\subset C^{m}(\,\overline{\Omega}\,)\cap
H^{\sigma+2q,1}(\Omega).
$$
(Последнее равенство становится ясным, если рассмотреть отдельно
случаи $m+n/2\geq\sigma+2q$ и $m+n/2<\sigma+2q$ и воспользоваться
пунктом г) теоремы 4.3.) Кроме того,
$$
\chi u\in H^{2q+n/2,\varphi,(2q)}(\Omega)=
H^{2q+n/2,\varphi}(\Omega)\subset C^{2q}(\,\overline{\Omega}\,)
$$
для любой функции $\chi\in C^{\infty}_{0}(\Omega)$, что влечет
включение $u\in C^{2q}(\Omega)$. Таким образом, выполняется
условие (6.8), т. е. $u$ --- классическое решение. Теорема 6.3
доказана.

\vspace{0.5cm}

\textbf{7. Корректность определения некоторых пространств.}
Покажем, что пространство (2.4) не зависит от использованного в
его определении параметра $\varepsilon$.

\textbf{Теорема 7.1.} \it Пусть $r\in\mathbb{N}$, $s\in E_{r}$ и
$\varphi\in\mathcal{M}$. Пространство
$$
H^{s,\varphi,(r)}(\Omega,\varepsilon):=
\bigl[\,H^{s-\varepsilon,\varphi,(r)}(\Omega),
H^{s+\varepsilon,\varphi,(r)}(\Omega)\,\bigr]_{1/2}
$$
не зависит с точностью до эквивалентности норм от параметра
$\varepsilon\in(0,1)$. \rm

\it\textbf{Доказательство.} \rm Предположим сначала, что $r=2q$
--- четное число. Рассмотрим какую-либо регулярную эллиптическую
краевую задачу (1.1), для которой пространства $N$ и $N^{+}$
тривиальны. (Например, задачу Дирихле для $A:=(1-\Delta^{q})$, где
$\Delta$ --- оператор Лапласа.) Согласно теореме 5.1 существует
топологический изоморфизм
$$
(A,B):\,H^{s,\varphi,(2q)}(\Omega,\varepsilon)\leftrightarrow
\mathcal{H}_{s,\varphi}(\Omega,\Gamma)\quad\mbox{при}\quad
0<\varepsilon<1.
$$
Отсюда немедленно следует теорема для четного $r=2q$.

Предположим далее, что число $r$ нечетное. В силу теоремы 4.2 для
любого числа $\sigma\in(-\infty,r+1/2)\setminus E_{r}$ существуют
изометрические изоморфизмы
$$
T_{r}:\,H^{\sigma,\varphi,(r)}(\Omega)\leftrightarrow
K_{\sigma,\varphi,(r)}(\Omega,\Gamma),
$$
$$
T_{r+1}:\,H^{\sigma,\varphi,(r+1)}(\Omega)\leftrightarrow
K_{\sigma,\varphi,(r+1)}(\Omega,\Gamma)=
K_{\sigma,\varphi,(r)}(\Omega,\Gamma)\times
H^{\sigma-r-1/2,\varphi}(\Gamma).
$$
Поэтому композиция отображений
$$
u\mapsto T_{r+1}\,u=:(u_{0},u_{1},\ldots,u_{r},u_{r+1})\mapsto
(T_{r}^{-1}(u_{0},u_{1},\ldots,u_{r}),u_{r+1}),\quad u\in
H^{\sigma,\varphi,(r+1)}(\Omega),
$$
определяет изометрический изоморфизм
$$
T:\,H^{\sigma,\varphi,(r+1)}(\Omega)\leftrightarrow
H^{\sigma,\varphi,(r)}(\Omega)\times
H^{\sigma-r-1/2,\varphi}(\Gamma).
$$
Возьмем здесь $\sigma=s\mp\varepsilon$, где $0<\varepsilon<1$, и
применим интерполяцию со степенным параметром $t^{1/2}$. Получим
топологический изоморфизм
$$
T:\,H^{s,\varphi,(r+1)}(\Omega,\varepsilon)\leftrightarrow
H^{s,\varphi,(r)}(\Omega,\varepsilon)\times
H^{s-r-1/2,\varphi}(\Gamma):=X(\varepsilon). \eqno(7.1)
$$
При этом используется предложение 3.2 и интерполяционное равенство
$$
\left[H^{s-\varepsilon-r-1/2,\varphi}(\Gamma),
H^{s+\varepsilon-r-1/2,\varphi}(\Gamma)\right]_{1/2}=
H^{s-r-1/2,\varphi}(\Gamma)\quad\mbox{с эквивалетностью норм},
$$
которое доказывается аналогично равенству (5.7). Теперь в силу
(7.1) имеем
$$
\bigl\|\,u\,\bigr\|_{H^{s,\varphi,(r)}(\Omega,\varepsilon)}=
\bigl\|\,(u,0)\,\bigr\|_{X(\varepsilon)}\asymp
\bigl\|\,T^{-1}(u,0)\,\bigr\|_{H^{s,\varphi,(r+1)}(\Omega,\varepsilon)}.
$$
Отсюда, поскольку параметр $r+1$ четный, вытекает, по доказанному,
что нормы в пространствах $H^{s,\varphi,(r)}(\Omega,\varepsilon)$,
где $0<\varepsilon<1$, эквивалентны. Значит, эти пространства
равны, поскольку множество $C^{\infty}(\,\overline{\Omega}\,)$
плотно в каждом из них согласно теореме 4.3 б). Теорема 7.1
доказана.

\end{document}